\pdfoutput=1
\documentclass[thmsa,notitlepage,leqno]{article}
\usepackage{amssymb}
\usepackage{amsfonts}
\usepackage[leqno]{amsmath}
\usepackage{sw20jart}
\usepackage{graphicx}
\usepackage[a4paper,margin=1in]{geometry}

\input{tcilatex}

\begin{document}

\author{Pablo Azcue\thanks{Departamento de Matematicas, Universidad Torcuato Di Tella.
Argentina.}, Esther Frostig\thanks{Department of Statistics, University of
Haifa. Israel.} and Nora Muler$^{\ast}$}
\title{Optimal strategies in a production-inventory control model}
\date{}
\maketitle

\begin{abstract}
We consider a production-inventory control model with finite capacity and two
different production rates, assuming that the cumulative process of customer
demand is given by a compound Poisson process. It is possible at any time to
switch over from the different production rates but it is mandatory to
switch-off when the inventory process reaches the storage maximum capacity. We
consider holding, production, shortage penalty and switching costs. This model
was introduced by Doshi, Van Der Duyn Schouten and Talman in 1978. Our aim is
to minimize the expected discounted cumulative costs up to infinity over all
admissible switching strategies. We show that the optimal cost functions for
the different production rates satisfy the corresponding
Hamilton-Jacobi-Bellman system of equations in a viscosity sense and prove a
verification theorem. The way in which the optimal cost functions solve the
different variational inequalities gives the switching regions of the optimal
strategy, hence it is stationary in the sense that depends only on the current
production rate and inventory level. We define the notion of finite band
strategies and derive, using scale functions, the formulas for the different
costs of the band strategies with one or two bands. We also show that there
are examples where the switching strategy presented by Doshi et al. is not the
optimal strategy.


\end{abstract}

\section{Introduction}

The classical production-inventory model considers a single machine that
produces a certain product. Finished products are stored and the storage
capacity can be finite or infinite. Moreover, the classical model assume a
constant production rate, customer demands arriving according to a Poisson
process and size demands distributed as i.i.d random variables. When the
stock on hand is less than the demand then either the excess of the demand
is lost or backlogged. In the first case the inventory level is always
positive, while in the latter it can be negative. The costs associated with
this model are holding cost and lost-sales cost. Higher production rates
yield fewer lost-sale cost but higher holding cost and viceversa. Thus,
there is a trade-off between holding and lost-sales costs. Therefore,
researchers have looked for the optimal strategy to minimize the expected
cost. One of the prominent strategy discussed in the literature is the two
regime switching policy. Under this policy, the production rate switches
from high to low rate when the inventory increases above a given level $%
y_{1} $; also, the production rate switches from low to high rate when the
inventory becomes smaller than a given level $y_{2},$ where $y_{2}<y_{1}$.

In the operations research literature, most articles have considered the
average cost per time unit assuming that the system is at steady state.
Gavish and Graves \cite{GG} and Gavish and Keilson \cite{GK} studied the
case where once the inventory level reaches a given threshold $y_{1},$
productions stops; and production resumes when the inventory level
down-crosses another threshold $y_{2}$, where $y_{2}<y_{1}$. In these two
papers, customers arrive according to Poisson process and backlogging is
permitted. In the first paper the demand is always for one item and the
machine produces one item per time unit, and in the second one the demand is
exponentially distributed. In both papers, the average cost per time unit is
obtained. De Kok, Tijms and Van Der Schouten \cite{kok84}, De Kok \cite%
{kok85} and De Kok \cite{kok87} \ studied an infinite capacity production
inventory system where demand occurs according to a compound Poisson process
and unsatisfied demand is backlogged. They considered two production rates $%
\sigma _{2}<\sigma _{1}$ where the production rate is switched to $\sigma
_{2}$ once the inventory level is above $y_{1}$ and it is switched back to $%
\sigma _{{1}}$ when the inventory level down-crosses $y_{2}$. In the first
paper, unsatisfied demand is backlogged and in the second one, unsatisfied
demand is lost. Performance measures that are considered under some
constrains on the switching and holding costs are: the average amount of
stock-out per unit time, the fraction of demand to be met directly from
stock on hand (in the backlog case) and the average amount of lost sales.
Doshi, Van Der Duyn and Talman \cite{Doshi} considered a finite capacity
production inventory model with lost sales and similar production rate
policy as in \cite{kok84}, \cite{kok85} and \cite{kok87}. They obtained the
steady-state distribution of the inventory level for this model and hence
the average cost per time unit.

More recently, Shi, Katehakis, Melamed and Xia \cite{shi} considered an
infinite capacity production-inventory model with compound Poisson demand,
lost sales and constant production rate. They obtained the expected
discounted cost and then the production rate that minimizes it. Barron,
Perry and Stadje \cite{barron} considered the model of Doshi et al. \cite%
{Doshi} under the assumption of Markov additive arrival process and
phase-type demand and obtained the expected discounted cost.

The optimal two-regime switching problem, also called starting-and-stopping
problem, has been studied extensively, in the diffusion setting and some
special profit functions, Brekke and Oksendal \cite{Brekke} apply a
verification approach for solving the variational inequality associated with
this impulse control problem. Pham and Vath \cite{PhamVath}, Hamad\`{e}ne and
Jeanblanc \cite{Hamadene}, and Bayraktar and Egami \cite{Bayraktar} between
others, studied various extensions of this model. Also in the diffusion
setting, Pham, Vath and Zhou \cite{Pham multiple} considered the case of
multiple-regime switching. Azcue and Muler \cite{AM Switching} studied a mixed
singular control/switching problem for multiple regimes in the compound
Poisson setting.

The rest of the paper is structured as follows. Section \ref{Section Model}
describes the model setup and some basic results are derived in Section
\ref{Section Basic Properties}. In Section \ref{Section HJB equations}, we
show that the optimal cost functions for the different production rates
satisfy the corresponding Hamilton-Jacobi-Bellman system of equations in a
viscosity sense and prove both characterization and verification results.
Moreover, we prove that there exists an optimal production-inventory strategy
and that it has a band structure. In Section
\ref{Section Finite Band Strategies}, we introduce the concept of finite band
strategies depending on the number of connected components of the non-action
regions; and in Section \ref{Section Value Functions} we use the scale
functions to find the formulas of the holding, shortage and switching cost
functions for the band strategies with one or two connected components.
Finally, in Section \ref{Section Examples}, we identify the optimal strategies
and the corresponding cost functions for a number of concrete examples with
exponentially distributed customer demands.

\section{Model\label{Section Model}}

In this paper we address a production-inventory control model with finite
storage capacity $b>0$ and two production rates: $\sigma_{1}$ and $\sigma_{2}$
such that $0<\sigma_{2}<\sigma_{1}$; this model was introduced by Doshi et al.
\cite{Doshi}. We say that the production is in phase $i=1,2$ when the
production rate is $\sigma_{i}$,\ whenever the inventory level reaches level
$b$, the production is stopped i.e. $\sigma_{0}=0$ at inventory level $b$. We
assume that the cumulative process of customer demand is given by the compound
Poisson process%

\[
\sum_{n=1}^{N_{t}}Y_{n}\text{,}%
\]
where $N_{t}$ is a Poisson process with rate of arrival $\lambda$ and the size
of the demand $Y_{n}$ are i.i.d positive random variables with distribution
$F$ and finite mean. Let us call $\tau_{n}$ the arrival of customer demand
$n$. \ We also assume that $l\leq0$ is the minimum level below which the
inventory is not allowed to decrease. If the inventory drops below $l$ the
part of the demand below $l$ is lost and production resumes at inventory level
$l$.

The following costs are considered:

\begin{itemize}
\item \textbf{Holding and production costs}. $h_{i}:[l,b)\rightarrow
\lbrack0,\infty)$ for $i=1,2$ correspond to the holding and production cost in
phase $i$ when the inventory level is $x\in\lbrack l,b)$. We assume that it is
bounded with finitely many discontinuities and Lipschitz between
discontinuities with Lipschitz constant $m_{h}$. $h_{0}(b)\geq0$ corresponds
to the holding cost at inventory level $b$.

\item \textbf{Shortage penalty costs}. $p:[0,\infty)\rightarrow\lbrack
0,\infty)$ corresponds to the penalty function cost when the amount $y$ of the
demand of a customer is lost. We assume that it is non-negative and
non-decreasing. Moreover,
\begin{equation}
\int_{0}^{\infty}p(y)dF(y)<\infty. \label{finite expectation penalty}%
\end{equation}

\item \textbf{Switching costs}. $K_{ij}$ corresponds to the fixing cost of
switching from phase $i$ to phase $j$ where $i$, $j=0,1,2$. Here we include
the costs of switch on ($K_{0i}$ where $i=1,2$) and the costs of switch off
($K_{i0}$ where $i=1,2$) the production process when the inventory reaches
level $b$. We add the following conditions on the switching costs in order to
penalize simultaneous changes of phases:
\begin{equation}%
\begin{array}
[c]{c}%
K_{0i}\leq K_{0j}+K_{ji}~\text{for }\{i,j\}=\{1,2\}\text{,}\\
\multicolumn{1}{l}{K_{12}+K_{21}>0.}%
\end{array}
\label{Condiciones K}%
\end{equation}

\end{itemize}

\begin{remark}
We assume here that it is possible at any time to switch over from phase $i$
to phase $j$ where $1\leq i$, $j\leq2$ but it is mandatory to switch off
(namely to go to phase $0$) when the inventory process reaches level $b$. On
top of that, the phase should be $1$ or $2$ (that corresponds to positive
production rate)\ whenever the inventory process is in the interval $\left[
l,b\right)  $. Moreover, if a demand arrives and the inventory level before
this arrival minus the demand of the customer is less than the backlog
$l\leq0,$ this demand is covered up to $l$ paying the corresponding penalty
cost of the part of the demand that has been lost given by function $p$ .
\end{remark}

Our aim is to minimize the expected discounted cumulative costs over all
possible production strategies. A production strategy can be defined as
$\pi=(T_{k},J_{k})_{k\geq1}$ where $T_{k}$ are the switching times from phase
$J_{k-1}$ to phase $J_{k}$ and $J_{k}\in\{0,1,2\}$. We call $T_{0}=0$ and
$J_{0}$ as the initial phase. In addition, we assume that $T_{1}<T_{2}<T_{3}$
$<\cdots$, and $J_{k}\neq J_{k-1}.$

Given a initial inventory level $x$, an initial phase $J_{0}=i$ and a
production strategy $\pi=(T_{k},J_{k})_{k\geq1},$ the controlled process is
defined recursively as $X_{T_{0}}^{\pi}=x$, $T_{0}=0$, and%

\begin{equation}
\text{ }X_{t}^{\pi}=X_{T_{k}}^{\pi}+\sigma_{J_{k}}\left(  t-T_{k}\right)
-\sum\nolimits_{n=N_{T_{k}}}^{N_{t}}\min\{Y_{n},X_{\tau_{n}^{-}}^{\pi
}-l\}\text{ for }t\in\lbrack T_{k},T_{k+1}). \label{controlled process}%
\end{equation}

Let us define the auxiliary inventory process,%

\[
\overset{\vee}{X}_{t}^{\pi}:=\text{ }X_{t}^{\pi}\text{ for }t\neq\tau
_{n}\text{ and }\overset{\vee}{X}_{\tau_{n}}^{\pi}=\text{ }X_{\tau_{n}^{-}%
}^{\pi}-Y_{n}\text{,}%
\]
so $X_{\tau_{n}}^{\pi}=l\vee\overset{\vee}{X}_{\tau_{n}}^{\pi}$, this
corresponds to the controlled process before it eventually resumes at
inventory level $l$.

Let us also define the controlled phase process
\begin{equation}
\mathcal{J}_{t}:=J_{k}\text{ for }t\in\lbrack T_{k},T_{k+1}). \label{Jt}%
\end{equation}

A production strategy $\pi=(T_{k},J_{k})_{k\geq0}\in\Pi_{x,i}$ starting at
phase $i$ and inventory level $x$ is \textit{admissible} if it is
$\mathcal{F}_{t}$-adapted, c\`{a}dl\`{a}g and satisfies,

\begin{itemize}
\item $T_{0}=0$ and $J_{0}=i.$

\item If the current inventory level is less than $b$, then the phase should
be either $1$ or $2.$ More precisely, if $X_{t}^{\pi}<b$ then $\mathcal{J}%
_{t^{-}}$ must be $1$ or $2$.

\item If at time $t,$ the phase process $\mathcal{J}_{t^{-}}=i$ with $i=1,2$
and the current inventory level $X_{t}^{\pi}$ level reaches $b,$ it is
mandatory to switch off the production. Hence, this time $t$ should coincide
with the next switching time $T_{k}$ for some $k$ and $\mathcal{J}_{t}%
=J_{k}=0$. Afterwards, $X_{t}^{\pi}=b$ for $t\in\lbrack T_{k},T_{k+1})$ , and
$T_{k+1}$ would be the time of the arrival of the next costumer demand and
$J_{k+1}$ would be either $1$ or $2$.
\end{itemize}

If the initial phase is $i\in\left\{  1,2\right\}  $, given an initial
inventory level $x\in\lbrack l,b)$, and an admissible production strategy
$\pi=(T_{k},J_{k})_{k\geq0}\in\Pi_{x,i}$, the associated cost function is
given by,%
\[%
\begin{array}
[c]{lll}%
V_{i}^{\pi}(x) & = & \mathbb{E}\left[  \int_{0}^{\infty}e^{-qt}h_{\mathcal{J}%
t}(X_{t}^{\pi})dt\right]  +\mathbb{E}\left[  \sum_{k=0}^{\infty}e^{-qT_{k+1}%
}K_{J_{k},J_{k+1}}\right] \\
&  & +\mathbb{E}\left[  \sum\nolimits_{n=1}^{\infty}e^{-q\tau_{n}}1_{\left\{
X_{\tau_{n}{}^{-}}^{\pi}-l<Y_{n}\right\}  }~p\left(  Y_{n}-X_{\tau_{n}{}^{-}%
}^{\pi}+l\right)  \right]  .
\end{array}
\]

We define the optimal cost functions for $i=1,2$ as%
\begin{equation}
V_{i}(x)=\inf_{\pi\in\Pi_{x,i}}V_{i}^{\pi}(x) \label{Definicion Vi}%
\end{equation}
for $x\in\lbrack l,b)$.

Given an initial inventory level $b$ and an admissible inventory strategy
$\pi=(T_{k},J_{k})_{k\geq0}\in\Pi_{b,0}$ the cost value of this strategy is
given by%

\[%
\begin{array}
[c]{lll}%
V_{0}^{\pi}(b) & = & \mathbb{E}\left[  \int_{0}^{\infty}e^{-qt}h_{\mathcal{J}%
t}(X_{t}^{\pi})dt\right]  +\mathbb{E}\left[  \sum_{k=0}^{\infty}e^{-qT_{k+1}%
}K_{J_{k},J_{k+1}}\right] \\
&  & +\mathbb{E}\left[  \sum\nolimits_{n=1}^{\infty}e^{-q\tau_{n}}1_{\left\{
X_{\tau_{n}{}^{-}}^{\pi}-Y_{n}<l\right\}  }~p\left(  Y_{n}-X_{\tau_{n}{}^{-}%
}^{\pi}+l\right)  \right]  .
\end{array}
\]
In this case, the optimal value for inventory level $b$ is given by%

\begin{equation}
V_{0}(b)=\inf_{\pi\in\Pi_{b,0}}V_{0}^{\pi}(b). \label{Definicion V0}%
\end{equation}

\section{Basic Properties\label{Section Basic Properties}}

In this section we study the existence and regularity of the optimal cost
functions. Let us start proving that they are well defined.

\begin{proposition}
\label{Proposicion Bien Definida} $V_{0}(b)$ is finite and the optimal cost
functions $V_{i}$ are bounded in $[l,b)$ for $i=1,2$. We call $\overline
{V}_{i}$ the positive upper bounds of the functions $V_{i}$ for $i=1,2.$
\end{proposition}

Proof.

Take $i\in\left\{  1,2\right\}  $, $x\in\lbrack l,b)$ and the admissible
production strategy $\pi=(T_{k},J_{k})_{k\geq1}\in\Pi_{x,i}$ that only switch
off from phase $i$ to phase $0$ when the current inventory level is $b$ and
remain in phase $i$ otherwise. Let us call
\begin{equation}
\overline{h}=\max\left\{  \sup_{x\in\lbrack l,b]}h_{1}(x),\sup_{x\in\lbrack
l,b]}h_{2}(x),h_{0}(b)\right\}  . \label{Definicion maximo h}%
\end{equation}
Then, we have%

\begin{equation}
\mathbb{E}\left[  \int_{0}^{\infty}e^{-qt}h_{\mathcal{J}t}(X_{t}^{\pi
})dt\right]  \leq\frac{\overline{h}}{q}. \label{cota h}%
\end{equation}
Moreover,%

\begin{equation}%
\begin{array}
[c]{l}%
\mathbb{E}\left[  \sum\nolimits_{n=1}^{N_{t}}e^{-q\tau_{n}}1_{\left\{
X_{\tau_{n}{}^{-}}^{\pi}<Y_{n}+l\right\}  }~p\left(  Y_{n}-X_{\tau_{n}{}^{-}%
}^{\pi}+l\right)  \right] \\
\leq\mathbb{E}\left[  \sum_{n=1}^{\infty}e^{-q\tau_{n}}p(Y_{n})\right]
=\mathbb{E}\left[  \sum_{n=1}^{\infty}e^{-q\tau_{n}}\right]  \mathbb{E}%
[p(Y_{1})],
\end{array}
\label{Cota p}%
\end{equation}
and so it is finite from (\ref{finite expectation penalty}). Finally,%

\begin{equation}%
\begin{array}
[c]{lll}%
\mathbb{E}\left[  \sum_{k=1}^{\infty}e^{-qT_{k}}K_{J_{k-1},J_{k}}\right]  &
\leq & \mathbb{E}\left[  \sum_{k=1}^{\infty}1_{\mathcal{J}_{T_{k}^{-}=i}%
}1_{\mathcal{J}_{T_{k}=0}}e^{-qT_{k}}K_{i,0}\right]  +\mathbb{E}\left[
\sum_{k=1}^{\infty}1_{\mathcal{J}_{T_{k}^{-}=0}}1_{\mathcal{J}_{T_{k}=i}%
}e^{-qT_{k}}K_{0,i}\right] \\
& \leq & K_{i,0}+\left(  K_{i,0}+K_{0,i}\right)  \mathbb{E}\left[  \sum
_{k=1}^{\infty}e^{-q\tau_{k}}\right] \\
& \leq & K_{i,0}+\left(  K_{i,0}+K_{0,i}\right)  \lambda/q
\end{array}
\label{Cota Switching}%
\end{equation}
so from (\ref{cota h}), (\ref{Cota p}) and (\ref{Cota Switching}), the
function $V_{i}$ is bounded in $[l,b)$. With a similar proof it can be shown
that $V_{0}$ is finite and so we have the result.$\blacksquare$

\begin{proposition}
The optimal cost functions $V_{i}$ are Lipschitz for $i=1,2$ in $[l,b)$.
\end{proposition}

Proof.

Given initial inventory level $x\in\lbrack l,b)$ and initial phase $i=1,2,$
take $\delta\in(0,b-x]$ and consider an admissible strategy $\pi_{x+\delta}%
\in\Pi_{x+\delta,i}$ such that $V_{i}^{\pi_{x+\delta}}(x+\delta)\leq
V_{i}(x+\delta)+\varepsilon,$ where $0<\varepsilon<\delta$. Let us now define
the admissible strategy $\pi_{x}\in\Pi_{x,i}$ as follows: stay in phase $i$
until the controlled inventory level $X_{t}^{\pi_{x}}$ reaches $x+\delta$ and
then follow $\pi_{x+\delta}\in\Pi_{x+\delta,i}$. Then, from
(\ref{Definicion maximo h}) and Proposition \ref{Proposicion Bien Definida},
we get%

\[%
\begin{array}
[c]{lll}%
V_{i}(x) & \leq & V_{i}^{\pi_{x}}(x)\\
& \leq & \int_{0}^{\frac{\delta}{\sigma_{i}}}e^{-qt}h_{\mathcal{J}t}%
(x+\sigma_{i}t)dt+\mathbb{P}[\tau_{1}>\frac{\delta}{\sigma_{i}}]e^{-q\frac
{\delta}{\sigma_{i}}}V_{i}^{\pi_{x+\delta}}(x+\delta)\\
&  & +\mathbb{P}\left[  \tau_{1}\leq\frac{\delta}{\sigma_{i}}\right]
\overline{V}_{i}\\
& \leq & \overline{h}\frac{\delta}{\sigma_{i}}+e^{-\left(  \lambda+q\right)
\frac{\delta}{\sigma_{i}}}\left(  V_{i}(x+\delta)+\varepsilon\right)
+(1-e^{-\lambda\frac{\delta}{\sigma_{i}}})\overline{V}_{i}.
\end{array}
\]
Hence, we have%
\[%
\begin{array}
[c]{lll}%
V_{i}(x)-V_{i}(x+\delta) & \leq & \overline{h}\frac{\delta}{\sigma_{i}%
}+e^{-\left(  \lambda+q\right)  \frac{\delta}{\sigma_{i}}}\left(
V_{i}(x+\delta)+\varepsilon\right)  -V_{i}(x+\delta)+(1-e^{-\lambda
\frac{\delta}{\sigma_{i}}})\overline{V}_{i}\\
& \leq & \overline{h}\frac{\delta}{\sigma_{i}}+\varepsilon+\lambda\frac
{\delta}{\sigma_{i}}\overline{V}_{i}.
\end{array}
\]
So, taking
\[
m_{i}^{1}:=\frac{\overline{h}}{\sigma_{i}}+1+\frac{\lambda}{\sigma_{i}%
}\overline{V}_{i},
\]
we obtain%

\begin{equation}
V_{i}(x)-V_{i}(x+\delta)\leq m_{i}^{1}\delta. \label{Lipschitz1}%
\end{equation}
Let us prove now that there exists $m_{i}^{2}>0$ such that,%

\begin{equation}
V_{i}(x+\delta)-V_{i}(x)\leq m_{i}^{2}\delta. \label{Lipschitz2}%
\end{equation}
We start showing that there exists $m$ such that,%

\begin{equation}
V_{i}(y)-V_{i}(l)\leq m\delta\label{Lipschitz2L}%
\end{equation}
for all $y\in\lbrack l,l+\delta]$. Given $\varepsilon>0$ and an initial
inventory level $l$, consider the strategy $\pi_{l}\in\Pi_{l,i}$ for $i=1,2$
such that $V_{i}^{\pi_{l}}(l)\leq V_{i}(l)+\varepsilon$ and call $X_{t}%
^{\pi_{l}}$ the associated process with initial inventory level $l$. Take also
a strategy $\pi_{b}\in\Pi_{b,0}$ such that $V_{0}^{\pi_{b}}(b)\leq
V_{0}(b)+\varepsilon$.

Let us define the admissible strategy $\pi_{y}\in\Pi_{y,i}$ for initial
inventory level $y\in\lbrack l,l+\delta]$ as:

\begin{itemize}
\item For $0\leq t\leq T,$ follow $\pi_{l}$ (and so the associated controlled
processes $X_{t}^{\pi_{y}}=X_{t}^{\pi_{l}}+(y-l)$ for $t<T$), where
\[
T:=\min\{t:X_{t}^{\pi_{y}}=b\text{ or }\overset{\vee}{X}_{t}^{\pi_{y}%
}-(y-l)=\overset{\vee}{X}_{t}^{\pi_{l}}<l\}.
\]

\item If $X_{T}^{\pi_{y}}=b$, follow $\pi_{b}$ for $t\geq T$.

\item If $\overset{\vee}{X}_{T}^{\pi_{y}}<l$ (and so $X_{T}^{\pi_{y}}%
=X_{T}^{\pi_{l}}=l$), follow $\pi_{l}$ for $t\geq T$.

\item If $l\leq X_{T}^{\pi_{y}}<y$ (and so $X_{T}^{\pi_{l}}=l$ and $X_{T}%
^{\pi_{y}}=\overset{\vee}{X}_{T}^{\pi_{y}}$), also follow the strategy
$\pi_{l}$ for $t\geq T$.
\end{itemize}

Given any stopping time $\tau$, let us define $\widehat{V}_{i}^{\pi_{y}%
}(y,\tau)$ as the expected discounted cost of the strategy before $\tau$ and
$\widetilde{V}_{i}^{\pi_{y}}(y,\tau)$ as the expected discounted cost of the
strategy after $\tau$. Thus,%

\[%
\begin{array}
[c]{l}%
V_{i}(y)-V_{i}(l)-\varepsilon\\%
\begin{array}
[c]{ll}%
\leq & V_{i}^{\pi_{y}}(y)-V_{i}^{\pi_{l}}(l)\\
\leq & \mathbb{E}\left[  \int_{0}^{T}e^{-qt}\left(  h_{\mathcal{J}_{t}}%
(X_{t}^{\pi_{l}}+(y-l))-h_{\mathcal{J}_{t}}(X_{t}^{\pi_{l}})\right)
dt\right]  +\\
& +\mathbb{E}\left[  1_{X_{T}^{\pi_{y}}=b}\left(  e^{-qT}\left(
K_{\mathcal{J}_{T},0}+V_{0}(b)+\varepsilon\right)  -\widetilde{V}_{i}^{\pi
_{l}}(l,T)\right)  \right] \\
& +\mathbb{E}\left[  1_{\overset{\vee}{X}_{T}^{\pi_{y}}<l}e^{-qT}\left(
V_{\mathcal{J}_{T}}(l)+\varepsilon+p(l-\overset{\vee}{X}_{T}^{\pi_{y}%
})-(V_{\mathcal{J}_{T}}(l)+p(l-\overset{\vee}{X}_{T}^{\pi_{y}}+y-l))\right)
\right] \\
& +\mathbb{E}\left[  1_{\left\{  l\leq X_{T}^{\pi_{y}}<y\right\}  }%
e^{-qT}\left(  V_{\mathcal{J}_{T}}^{\pi_{X_{T}^{\pi_{y}}}}(X_{T}^{\pi_{y}%
})-V_{\mathcal{J}_{T}}^{\pi_{l}}(l)+2\varepsilon\right)  \right]  .
\end{array}
\end{array}
\]

Let $n_{D}$ be the sum of the numbers of discontinuities of $h_{1}$ and
$h_{2}$. Note that between two customer demands, the inventory level
$X_{t}^{\pi_{y}}$ goes through at most $n_{D}$ points of discontinuities of
$h_{\mathcal{J}_{t}}$. Hence, calling $\tau_{0}=0$, we have%

\begin{equation}
\mathbb{E}\left[  \int_{0}^{T}e^{-qt}\left(  h_{\mathcal{J}t}(X_{t}^{\pi_{l}%
}+(y-l))-h_{\mathcal{J}t}(X_{t}^{\pi_{l}})\right)  dt\right]  \leq\left(
\frac{m_{h}}{q}+n_{D}\frac{\overline{h}}{\sigma_{2}}\left(  1+\frac{\lambda
}{q}\right)  \right)  \delta. \label{LL0}%
\end{equation}
Let us call $\widetilde{T}:=\inf\left\{  t:X_{t}^{\pi_{l}}=b\right\}  $ and
$\widetilde{\tau}$ the time of the first customer demand after $T;$ we have
that $\mathbb{P}\left[  \widetilde{T}>\widetilde{\tau}\right]  \leq
1-e^{-\lambda\frac{\delta}{\sigma_{2}}}$ and so%

\begin{equation}%
\begin{array}
[c]{l}%
\mathbb{E}\left[  1_{X_{T}^{\pi_{y}}=b}\left(  e^{-qT}\left(  K_{\mathcal{J}%
_{T},0}+V_{0}(b)\right)  -\widetilde{V}_{i}^{\pi_{l}}(l,T)\right)  \right] \\%
\begin{array}
[c]{ll}%
\leq & \left(  1-e^{-\lambda\frac{\delta}{\sigma_{2}}}\right)  \left(
V_{0}(b)+\left(  K_{1,0}\vee K_{2,0}\right)  \right) \\
& +\mathbb{E}\left[  1_{X_{T}^{\pi_{y}}=b}1_{\widetilde{T}<\widetilde{\tau}%
}e^{-qT}\left(  K_{\mathcal{J}_{T},0}+V_{0}(b)\right)  -\widetilde{V}_{i}%
^{\pi_{l}}(l,T)\right] \\
\leq & \lambda\frac{\delta}{\sigma_{2}}\left(  V_{0}(b)+\max\left\{
K_{1,0},K_{2,0}\right\}  \right) \\
& +\mathbb{E}\left[  1_{X_{T}^{\pi_{y}}=b}1_{\widetilde{T}<\widetilde{\tau}%
}e^{-qT}\left(  K_{\mathcal{J}_{T},0}+V_{0}(b)\right)  -\widetilde{V}_{i}%
^{\pi_{l}}(l,T)\right]  .
\end{array}
\end{array}
\label{LL1}%
\end{equation}

Let $\Delta$ be the length of time after $T$ in which the process $X_{t}%
^{\pi_{l}}$ reaches $b$ in the event of no arrivals of demands. In this case,
we have%
\[
X_{T+\Delta}^{\pi_{l}}=b-(y-l)+\int_{T}^{T+\Delta}e^{-qs}\sigma_{\mathcal{J}%
_{T}}ds=b
\]
and so $\frac{\delta}{\sigma_{1}}\leq\Delta\leq\frac{\delta}{\sigma_{2}}$.
Hence, from (\ref{Condiciones K}),%

\[%
\begin{array}
[c]{l}%
\mathbb{E}\left[  1_{\left\{  X_{T}^{\pi_{y}}=b\right\}  }1_{\left\{
\widetilde{T}<\widetilde{\tau}\right\}  }\widetilde{V}_{i}^{\pi_{l}%
}(l,T)\right] \\%
\begin{array}
[c]{ll}%
\geq & \mathbb{P}\left[  \text{no demands in }t\in\left[  T,T+\frac{\delta
}{\sigma_{2}}\right]  \right]  \text{ }\mathbb{E}\left[  e^{-q\left(
T+\Delta\right)  }\left(  K_{\mathcal{J}_{T},0}+V_{0}(b)\right)  \right] \\
\geq & e^{-\left(  q+\lambda\right)  \frac{\delta}{\sigma_{2}}}\mathbb{E}%
\left[  e^{-qT}\left(  K_{\mathcal{J}_{T},0}+V_{0}(b)\right)  \right]  .
\end{array}
\end{array}
\]
Therefore,%

\begin{equation}%
\begin{array}
[c]{l}%
\mathbb{E}\left[  1_{\left\{  X_{T}^{\pi_{y}}=b\right\}  }1_{\left\{
\widetilde{T}<\tau_{1}\right\}  }e^{-qT}\left(  K_{\mathcal{J}_{T},0}%
+V_{0}(b)\right)  -\widetilde{V}_{i}^{\pi_{l}}(l,T)\right] \\%
\begin{array}
[c]{ll}%
\leq & \left(  1-e^{-\left(  q+\lambda\right)  \frac{\delta}{\sigma_{2}}%
}\right)  \left(  K_{1,0}\vee K_{2,0}+V_{0}(b)\right) \\
\leq & \frac{q+\lambda}{\sigma_{2}}\left(  K_{1,0}\vee K_{2,0}+V_{0}%
(b)\right)  \delta.
\end{array}
\end{array}
\label{LL2}%
\end{equation}
Since the penalty function $p$ is non-decreasing, we also have,%
\begin{equation}
\mathbb{E}\left[  1_{\left\{  \overset{\vee}{X}_{T}^{\pi_{y}}<l\right\}
}e^{-qT}\left(  p(l-\overset{\vee}{X}_{T}^{\pi_{y}})-p(l-\overset{\vee}{X}%
_{T}^{\pi_{y}}+y-l)\right)  \right]  \leq0. \label{LL3nueva}%
\end{equation}

Finally, since the event $l\leq X_{T}^{\pi_{y}}<y$ coincides with the arrival
of a customer demand,%

\begin{equation}%
\begin{array}
[c]{l}%
\mathbb{E}\left[  1_{\left\{  l\leq X_{T}^{\pi_{y}}<y\right\}  }e^{-qT}\left(
V_{\mathcal{J}_{T}}^{\pi_{X_{T}^{\pi_{y}}}}(X_{T}^{\pi_{y}})-V_{\mathcal{J}%
_{T}}(l)\right)  \right] \\%
\begin{array}
[c]{ll}%
= & \mathbb{E}\left[  1_{\left\{  l\leq X_{T}^{\pi_{y}}<y\right\}
}1_{\left\{  T=\tau_{k}\text{ for some }k\right\}  }e^{-qT}\left(
V_{\mathcal{J}_{T}}^{\pi_{X_{T}^{\pi_{y}}}}(X_{T}^{\pi_{y}})-V_{\mathcal{J}%
_{T}}(l)\right)  \right] \\
\leq & \mathbb{E}\left[  e^{-q\tau_{1}}\max_{z\in\lbrack l,y]}\left(
V_{\mathcal{J}_{T}}^{\pi_{z}}(z)-V_{\mathcal{J}_{T}}(l)\right)  \right] \\
\leq & \frac{\lambda}{q+\lambda}\max_{z\in\lbrack l,l+\delta]}\left(
V_{\mathcal{J}_{T}}^{\pi_{z}}(z)-V_{\mathcal{J}_{T}}(l)\right)  .
\end{array}
\end{array}
\label{LL4}%
\end{equation}

Hence, from (\ref{LL1}), (\ref{LL2}), (\ref{LL3nueva}) and (\ref{LL4}), there
exists $\overline{m}$ large enough such that%

\[
\frac{q}{q+\lambda}\max_{z\in\lbrack l,l+\delta]}\left(  V_{i}^{\pi_{y}%
}(z)-V_{i}^{\pi_{l}}(l)\right)  \leq\overline{m}\delta.
\]
So, we obtain (\ref{Lipschitz2L}) with $m=\overline{m}\left(  q+\lambda
\right)  /q$. \bigskip The argument to show (\ref{Lipschitz2}) is analogous.
$\blacksquare$

\section{Hamilton Jacobi Bellman equations \label{Section HJB equations}}

From the definitions (\ref{Definicion Vi}) and (\ref{Definicion V0}), we can
obtain recursive equations relating the optimal cost $V_{0}(b)\ $and the
optimal cost functions $V_{i}$ for $i=1,2$; these recursive equations will be
used to find the Hamilton-Jacobi-Bellman equations of the optimization problem.

It follows immediately from (\ref{Definicion V0}) that%

\begin{equation}%
\begin{array}
[c]{lll}%
V_{0}(b) & = & \mathbb{E}\left[  \left(  \int_{0}^{\tau_{1}}e^{-qs}%
h_{0}(b)ds+1_{\left\{  Y_{1}\leq b-l\right\}  }e^{-q\tau_{1}}\overline
{V}(b-Y_{1})\right)  \right] \\
&  & +\mathbb{E}\left[  \left(  1_{\left\{  Y_{1}>b-l\right\}  }e^{-q\tau_{1}%
}\left(  p(Y_{1}-b+l)+\overline{V}(l)\right)  \right)  \right] \\
& = & \frac{1}{q+\lambda}h_{0}(b)+\frac{\lambda}{q+\lambda}%
{\textstyle\int_{0}^{b-l}}
\overline{V}(b-\alpha)dF(\alpha)\\
&  & +\frac{\lambda}{q+\lambda}\left(
{\textstyle\int_{b-l}^{\infty}}
p(\alpha-b+l)dF(\alpha)+\overline{V}(l)(1-F(b-l))\right)  ,
\end{array}
\label{DPP0}%
\end{equation}
where%
\begin{equation}
\overline{V}(x)=\min\{K_{01}+V_{1}(x),K_{02}+V_{2}(x)\}.
\label{Definicion Vbarra}%
\end{equation}
For $x\in\lbrack l,b)$, let us define
\[
t_{x}^{i}:=\min\{t:x+\sigma_{i}t=b\}=\frac{b-x}{\sigma_{i}}.
\]
Take $\left\{  i,j\right\}  =\left\{  1,2\right\}  $ and consider any stopping
time $T_{1}\geq0$ and $0<h<t_{x}^{i}$. Define $\tau=\tau_{1}\wedge T_{1}\wedge
h$ and
\[%
\begin{array}
[c]{lll}%
P_{i}(x,T_{1},h) & = & \mathbb{E}\left[  1_{\left\{  \tau=h<\tau_{1}\wedge
T_{1}\right\}  }\left(  \int_{0}^{h}e^{-qs}h_{i}(x+\sigma_{i}s)ds+V_{i}%
(x+\sigma_{i}h)e^{-qh}\right)  \right]  +\\
&  & +\mathbb{E}\left[  1_{\left\{  \tau=\tau_{1}<T_{1}\wedge h\right\}
}\left(  \int_{0}^{\tau_{1}}e^{-qs}h_{i}(x+\sigma_{i}s)ds\right)  \right] \\
&  & +\mathbb{E}\left[  1_{\left\{  \tau=\tau_{1}<T_{1}\wedge h\right\}
}1_{\left\{  Y_{1}\leq x+\sigma_{i}\tau_{1}-l\right\}  }e^{-q\tau_{1}}%
V_{i}(x+\sigma_{i}\tau_{1}-Y_{1})\right] \\
&  & +\mathbb{E}\left[  1_{\left\{  \tau=\tau_{1}<T_{1}\wedge h\right\}
}1_{\left\{  Y_{1}^{i}>x+\sigma_{i}\tau_{1}-l\right\}  }e^{-q\tau_{1}}\left(
p\left(  Y_{1}-\left(  x+\sigma_{i}\tau_{1}^{i}-l\right)  \right)
+V_{i}(l)\right)  \right] \\
&  & +\mathbb{E}\left[  1_{\left\{  \tau=T_{1}<\tau_{1}\wedge h\right\}
}\left(  \int_{0}^{T_{1}}e^{-qs}h_{i}(x+\sigma_{i}s)ds+\left(  V_{j}%
(x+\sigma_{i}T_{1})+K_{i,j}\right)  e^{-qT_{1}}\right)  \right]  .
\end{array}
\]
We obtain the following recursive equations%

\begin{equation}
V_{i}(x)=\inf_{T_{1}\geq0}P_{i}(x,T_{1},h). \label{DPPi}%
\end{equation}

Let us define the operators,%

\begin{equation}%
\begin{array}
[c]{lll}%
\mathcal{L}_{i}(V_{i})(x) & := & \sigma_{i}V_{i}^{\prime}(x)-(\lambda
+q)V_{i}(x)+\lambda%
{\textstyle\int_{0}^{x-l}}
V_{i}(x-\alpha)dF(\alpha)+\lambda%
{\textstyle\int_{x-l}^{\infty}}
p(\alpha-x+l)dF(\alpha)\\
&  & +\lambda V_{i}(l)(1-F(x-l))+h_{i}(x)
\end{array}
\label{Operador Li}%
\end{equation}
for $i=1,2$. Then, the Hamilton-Jacobi-Bellman equations for $V_{i}$, are%

\begin{equation}
\min\{\mathcal{L}_{i}(V_{i})(x),V_{j}(x)+K_{ij}-V_{i}(x)\}=0, \label{HJBi}%
\end{equation}
for $x\in\lbrack l,b)$, $\left\{  i,j\right\}  =\left\{  1,2\right\}  $. Also,
defining
\begin{equation}%
\begin{array}
[c]{ccc}%
\mathcal{L}_{0}(V_{0})(b) & := & -(q+\lambda)V_{0}(b)+\lambda\left(
{\textstyle\int_{0}^{b-l}}
\overline{V}(b-\alpha)dF(\alpha)+%
{\textstyle\int_{b-l}^{\infty}}
p(\alpha-b+l)dF(\alpha)\right) \\
&  & +\lambda\overline{V}(l)(1-F(b-l)+h_{0}(b),
\end{array}
\label{Operator L0}%
\end{equation}
we obtain from (\ref{DPP0}), that
\begin{equation}
\mathcal{L}_{0}(V_{0})(b)=0\text{.} \label{HJB0}%
\end{equation}

\begin{definition}
\label{Viscosity}A function $\underline{u}_{i}:$ $[l,b]\rightarrow\mathbf{R}%
$\ is a \textit{viscosity subsolution}\ of (\ref{HJBi}) at $x\in\lbrack l,b)$
for $\{i,j\}=$\ $\left\{  1,2\right\}  $ if it is Lipschitz and any
continuously differentiable function $\psi_{i}:[l,b]\rightarrow\mathbf{R}%
\ $with $\psi_{i}(x)=\underline{u}_{i}(x)$ such that $\underline{u}_{i}%
-\psi_{i}$\ reaches the minimum at $x$\ satisfies
\[
\min\{\mathcal{L}_{i}(\psi_{i})(x),V_{j}(x)+K_{ij}-\underline{u}_{i}%
(x)\}\leq0.
\]
A function $\overline{u}_{i}:[l,b]\rightarrow\mathbf{R}$\ is a
\textit{viscosity supersolution} of (\ref{HJBi}) at $x\in\lbrack l,b)$ for
$\{i,j\}=$\ $\left\{  1,2\right\}  $\ if it is Lipschitz and any continuously
differentiable function $\varphi_{i}:[l,b]\rightarrow\mathbf{R}\ $with
$\varphi_{i}(x)=\overline{u}_{i}(x)$ and such that $\overline{u}_{i}%
-\varphi_{i}$\ reaches the maximum at $x$\ satisfies
\[
\min\{\mathcal{L}_{i}(\varphi_{i})(x),V_{j}(x)+K_{ij}-\overline{u}%
_{i}(x)\}\geq0.\
\]
The functions $\psi_{i}$ and $\varphi_{i}$ are called test-functions for
subsolution and supersolution respectively. If a function $u_{i}$ is both a
subsolution and a supersolution at $x$ it is called a \textit{viscosity
solution} of (\ref{HJBi}) at $x$.
\end{definition}

Crandall and Lions \cite{CL} introduced the concept of viscosity solutions for
first-order Hamilton-Jacobi equations. It is the standard tool for studying
HJB equations, see for instance Fleming and Soner \cite{FS}.

\begin{proposition}
\label{V is viscosity solution} The optimal cost functions $V_{i}$ satisfy
(\ref{HJBi}) in a viscosity sense, for $x\in\lbrack l,b)$ and $i=1,2$.
\end{proposition}

Proof.

Consider $\{i,j\}=\{1,2\}$, taking $x\in\lbrack l,b)$ and $T_{1}=0$ in
(\ref{DPPi}), it follows that $V_{j}(x)+K_{ij}-V_{i}(x)\geq0$. Take now
$0<h<T_{1}$ and $h<t_{x}^{i}$. \ From (\ref{DPPi}) and using that $\varphi
_{i}$ is a test-functions for supersolution%

\[%
\begin{array}
[c]{lll}%
\varphi_{i}(x)=V_{i}(x) & \leq & \mathbb{E}\left[  1_{h<\tau_{1}}\left(
\int_{0}^{h}e^{-qs}h_{i}(x+\sigma_{i}s)ds+\varphi_{i}(x+\sigma_{i}%
h)e^{-qh}\right)  \right] \\
&  & +\mathbb{E}\left[  1_{\tau_{1}<h}\left(  \int_{0}^{\tau_{1}}e^{-qs}%
h_{i}(x+\sigma_{i}s)ds\right)  \right] \\
&  & +\mathbb{E}\left[  1_{\tau_{1}<h}1_{\left\{  Y_{1}\leq x+\sigma_{i}%
\tau_{1}-l\right\}  }e^{-q\tau_{1}}\varphi_{i}(x+\sigma_{i}\tau_{1}%
-Y_{1})\right] \\
&  & +\mathbb{E}\left[  1_{\tau_{1}<h}1_{\left\{  Y_{1}^{i}>x+\sigma_{i}%
\tau_{1}-l\right\}  }e^{-q\tau_{1}}\left(  p\left(  Y_{1}-\left(  x+\sigma
_{i}\tau_{1}^{i}-l\right)  \right)  +\varphi_{i}(l)\right)  \right]  .
\end{array}
\]

Hence,%

\[%
\begin{array}
[c]{lll}%
0 & \leq & \mathbb{E}\left[  1_{h<\tau_{1}}\left(  \int_{0}^{h}e^{-qs}%
h_{i}(x+\sigma_{i}s)ds+\varphi_{i}(x+\sigma_{i}h)e^{-qh}\right)  \right]
-\varphi_{i}(x)\\
&  & +\mathbb{E}\left[  1_{\tau_{1}<h}\left(  \int_{0}^{\tau_{1}}e^{-qs}%
h_{i}(x+\sigma_{i}s)ds\right)  \right] \\
&  & +\mathbb{E}\left[  1_{\tau_{1}<h}1_{\left\{  Y_{1}\leq x+\sigma_{i}%
\tau_{1}-l\right\}  }e^{-q\tau_{1}}\varphi_{i}(x+\sigma_{i}\tau_{1}%
-Y_{1})\right] \\
&  & +\mathbb{E}\left[  1_{\tau_{1}<h}1_{\left\{  Y_{1}^{i}>x+\sigma_{i}%
\tau_{1}-l\right\}  }e^{-q\tau_{1}}\left(  p\left(  Y_{1}-\left(  x+\sigma
_{i}\tau_{1}^{i}-l\right)  \right)  +\varphi_{i}(l)\right)  \right]  ,
\end{array}
\]
and so, dividing by $h$ and taking $h\rightarrow0^{+},$ we obtain
$\mathcal{L}_{i}(\varphi_{i})(x)\geq0$. Hence $V_{i}$ is a viscosity
supersolution of (\ref{HJBi}) at $x$.

Let us prove now that $V_{i}$ is a viscosity subsolution of (\ref{HJBi}) at
any $x\in(l,b)$ for $i=1,2$. It is enough to consider the case $V_{j}%
(x)+K_{ij}-V_{i}(x)>0$. Arguing by contradiction, we assume that $V_{i}$ is
not a subsolution of (\ref{HJBi}) at $x$. We can find, as in Proposition 3.1
in Azcue and Muler \cite{AM Libro}, values $\varepsilon>0$, $h\in
(0,(x-l)/2\wedge(b-x)/2)$ and a continuously differentiable function $\psi
_{i}\geq V_{i}$ in $[0,y+h]$ with $\psi_{i}(x)=V_{i}(x)$ such that%

\[%
\begin{array}
[c]{l}%
V_{j}(y)+K_{ij}-V_{i}(y)\geq0\text{ for }y\in\lbrack l,b),\\
\mathcal{L}_{i}(\psi_{i})(y)\geq2q\varepsilon\text{ for }y\in\lbrack
x-h,x+h],\\
V_{i}(y)\geq\psi_{i}(y)+3\varepsilon\text{ for }y\in\lbrack l,x-h]\cup\{x+h\},
\end{array}
\]
and also
\[
V_{j}(y)+K_{ij}-V_{i}(y)>2\varepsilon\text{ for }y\in\lbrack x-h,x+h].
\]
Since $\psi_{i}$ is continuously differentiable we can find a positive
constant $C$ such that $\mathcal{L}_{i}(\psi_{i})(y)\leq C$ for all
$y\in\lbrack l,b)$.

Let us take any admissible production strategy $\pi=(T_{k},J_{k})_{k\geq1}%
\in\Pi_{x,i}$, consider the uncontrolled inventory process $X_{t}^{\pi}$ defined in
(\ref{controlled process}), and define the stopping times
\[
\overline{\tau}=\inf\{t>0:X_{t}^{\pi}\geq x+h\},\underline{\tau}=\inf
\{t>0:X_{t}^{\pi}\leq x-h\},
\]
and $\tau^{\ast}=T_{1}\wedge\underline{\tau}\wedge\overline{\tau}$. We get
that if $\tau^{\ast}=T_{1}$ and $X_{\tau^{\ast}}^{\pi}\in(x-h,x+h)$ then%
\[
V_{j}(X_{T_{1}}^{\pi})+K_{ij}\geq V_{i}(X_{T_{1}}^{\pi})+2\varepsilon\geq
\psi_{i}(X_{T_{1}}^{\pi})+2\varepsilon,
\]
and in the case that either $\tau^{\ast}<T_{1}\ $or $\tau^{\ast}=T_{1}$ and
$X_{\tau^{\ast}}^{\pi}\notin(x-h,x+h)$ we have that
\[
V_{i}(X_{\tau^{\ast}}^{\pi})\geq\psi_{i}(X_{\tau^{\ast}}^{\pi})+2\varepsilon.
\]
Since the function $e^{-qt}\psi_{i}(x)$ is continuously differentiable, using
the expression (\ref{controlled process}) and the change of variables formula
for finite variation processes (see Protter \cite{Protter}), we can write%
\begin{equation}%
\begin{array}
[c]{ll}%
\psi_{i}(X_{\tau^{\ast}}^{\pi})e^{-q\tau^{\ast}}-\psi_{i}(x)= & \int
\nolimits_{0}^{\tau^{\ast}}\psi_{i}^{\prime}(X_{s^{-}}^{\pi})e^{-qs}\sigma
_{i}ds-q\int\nolimits_{0}^{\tau^{\ast}}\psi_{i}(X_{s^{-}}^{\pi})e^{-qs}ds\\
& +\sum\limits_{X_{s^{-}}\neq X_{s},s\leq\tau^{\ast}}\left(  \psi_{i}%
(X_{s^{-}}^{\pi}-\Delta X_{s})-\psi_{i}(X_{s^{-}}^{\pi})\right)  e^{-qs},
\end{array}
\label{SPS1}%
\end{equation}
where $\Delta X_{s}=X_{s}-X_{s^{-}}$.

On the other hand, $X_{s}\neq X_{s^{-}}$ only at the arrival of a demand, so%

\begin{equation}%
\begin{array}
[c]{cl}%
M_{t}= &
{\textstyle\sum\limits_{X_{s^{-}}\neq X_{s},s\leq t}}
\left(  \psi_{i}(X_{s^{-}}^{\pi}-\Delta X_{s})-\psi_{i}(X_{s^{-}}^{\pi
})\right)  e^{-qs}\\
& -\lambda\int\nolimits_{0}^{t}e^{-qs}\int\nolimits_{0}^{\infty}\left(
\psi_{i}(X_{s^{-}}^{\pi}-\alpha)-\psi_{i}(X_{s^{-}}^{\pi})\right)
dF(\alpha)ds
\end{array}
\label{SPSmartingala}%
\end{equation}
is a martingale with zero-expectation, here we extend the definition of
$\psi_{i}$ for $y<l$ as $\psi_{i}(y)=$ $p(l-y)+\psi_{i}(l).$ Therefore,\- we
can combine (\ref{SPS1}) and (\ref{SPSmartingala}) to obtain
\begin{equation}
\psi_{i}(X_{\tau^{\ast}}^{\pi})e^{-q\tau^{\ast}}-\psi_{i}(x)=%
{\textstyle\int\nolimits_{0}^{\tau^{\ast}}}
\mathcal{L}_{i}(\psi_{i})(X_{s^{-}}^{\pi})e^{-qs}ds+M_{\tau^{\ast}}-%
{\textstyle\int\nolimits_{0}^{\tau^{\ast}}}
h_{i}(X_{s^{-}}^{\pi})e^{-qs}ds. \label{nueva1}%
\end{equation}
In the case that $\tau^{\ast}=T_{1}$ and $X_{T_{1}}^{\pi}\in(x-h,x+h)$, we
have from (\ref{nueva1}) that%

\[
\left(  V_{j}(X_{T_{1}}^{\pi})+K_{ij}\right)  e^{-q\tau^{\ast}}+%
{\textstyle\int\nolimits_{0}^{\tau^{\ast}}}
h_{i}(X_{s^{-}}^{\pi})e^{-qs}ds\geq V_{i}(x)+2\varepsilon+M_{\tau^{\ast}}.
\]
In the case that either $\tau^{\ast}<T_{1}\ $or $\tau^{\ast}=T_{1}$ and
$X_{\tau^{\ast}}^{\pi}\notin(x-h,x+h)$, we get%

\begin{equation}%
\begin{array}
[c]{lll}%
(V_{i}(X_{\tau^{\ast}}^{\pi})-2\varepsilon)e^{-q\tau^{\ast}}-V_{i}(x) & \geq &
\psi_{i}(X_{\tau^{\ast}}^{\pi})e^{-q\tau^{\ast}}-\psi_{i}(x)\\
& = &
{\textstyle\int\nolimits_{0}^{\tau^{\ast}}}
\mathcal{L}_{i}(\psi_{i})(X_{s^{-}}^{\pi})e^{-qs}ds+M_{\tau^{\ast}}-%
{\textstyle\int\nolimits_{0}^{\tau^{\ast}}}
h_{i}(X_{s^{-}}^{\pi})e^{-qs}ds.\\
& \geq &
{\textstyle\int\nolimits_{0}^{\tau^{\ast}}}
2q\varepsilon e^{-qs}ds+M_{\tau^{\ast}}-%
{\textstyle\int\nolimits_{0}^{\tau^{\ast}}}
h_{i}(X_{s^{-}}^{\pi})e^{-qs}ds.\\
& = & 2\varepsilon(1-e^{-q\tau^{\ast}})+M_{\tau^{\ast}}-%
{\textstyle\int\nolimits_{0}^{\tau^{\ast}}}
h_{i}(X_{s^{-}}^{\pi})e^{-qs}ds.
\end{array}
\label{nueva 2}%
\end{equation}
and so, by (\ref{nueva1}) and (\ref{nueva 2}),
\[
e^{-q\tau^{\ast}}V_{i}(X_{\tau^{\ast}}^{\pi})+%
{\textstyle\int\nolimits_{0}^{\tau^{\ast}}}
h_{i}(X_{s^{-}}^{\pi})e^{-qs}ds\geq V_{i}(x)+\varepsilon+M_{\tau^{\ast}}^{1}.
\]
Finally, we obtain that $V_{i}^{\pi}(x)\geq V_{i}(x)+2\varepsilon$, and this
contradicts the definition of $V_{i}$. $\blacksquare$

In the following proposition, we prove that the optimal cost functions are the
largest viscosity supersolutions of their corresponding HJB equations with
suitable boundary conditions.

\begin{proposition}
\label{Largest Viscosity Supersolution} Fix $x\in\lbrack l,b)$ and $j=1,2$ or
$x=b$ and $j=0$. Let $\overline{u}_{1}$ and $\overline{u}_{2}$\ be
non-negative viscosity supersolution of the corresponding HJB equation
(\ref{HJBi}) in $[l,b)$ and consider any admissible strategy $\pi=(T_{k}%
,J_{k})_{k\geq0}\in\Pi_{x,j}$. Defining%
\[
\overline{u}(x)=\min\{K_{01}+\overline{u}_{1}(x),K_{02}+\overline{u}_{2}(x)\}
\]
and since $\mathcal{L}_{0}(\overline{u}_{0})(b)=0,$
\end{proposition}

\[
\overline{u}_{0}(b)=\frac{\lambda}{q+\lambda}\left(
{\textstyle\int_{0}^{b-l}}
\overline{u}(b-\alpha)dF(\alpha)+%
{\textstyle\int_{b-l}^{\infty}}
\left(  p(\alpha-b+l)+\overline{u}(l)\right)  dF(\alpha)\right)  +\frac
{h_{0}(b)}{q+\lambda}.
\]
If we assume that%

\[
\overline{u}_{1}(b)\leq\overline{u}_{0}(b)+K_{10},\text{ }\overline{u}%
_{2}(b)\leq\overline{u}_{0}(b)+K_{20},
\]
then $\overline{u}_{j}(x)\leq V_{j}^{\pi}(x)$ for $j=1,2$ and $\overline
{u}_{0}(b)\leq V_{0}^{\pi}(b).$

Proof.

Consider $\pi\in\Pi_{x,j}$. Let us extend $\overline{u}_{1}$ and $\overline
{u}_{2}$\ as $\overline{u}_{i}(x)=\overline{u}_{i}(l)$ and $\overline{u}%
_{0}(x)=\overline{u}_{0}(l)$ for $x<l$. Consider the controlled risk process
$X_{t}^{\pi}$ starting at $x$ and the function $\mathcal{J}_{t}\ $defined in
(\ref{Jt}). Since $\overline{u}_{i}$ is Lipschitz for $i=1,2$, we obtain that
the function $t\rightarrow e^{-qt}~\overline{u}_{\mathcal{J}_{t}}(X_{t}^{\pi
})$ is absolutely continuous in between the stopping times $\left\{
0\right\}  \cup\left\{  \tau_{n}:n\geq1\right\}  \cup\left\{  T_{k}%
:k\geq1\right\}  $. So, taking
\[
m_{t}:=\max\{k:T_{k}\leq t\}\text{,}%
\]
we have
\begin{equation}%
\begin{array}
[c]{l}%
\overline{u}_{\mathcal{J}_{t}}(X_{t}^{\pi})e^{-qt}-\overline{u}_{j}(x)\\%
\begin{array}
[c]{ll}%
= & \sum_{k=0}^{m_{t}-1}\left(  \overline{u}_{_{J_{k+1}}}(X_{T_{k+1}}^{\pi
})e^{-qT_{k+1}}-\overline{u}_{_{J_{k}}}(X_{T_{k}}^{\pi})e^{-qT_{k}}\right)
+(\overline{u}_{J_{m_{t}}}(X_{t}^{\pi})e^{-qt}-\overline{u}_{_{J_{m_{t}}}%
}(X_{T_{m_{t}}}^{\pi})e^{-qT_{m_{t}}}).
\end{array}
\end{array}
\label{Expresion Inicial}%
\end{equation}

Let us define
\begin{equation}%
\begin{array}
[c]{lll}%
M^{i}(z_{0},t_{0},t) & = & \overline{u}_{i}(Z_{t}^{i})e^{-qt}-\overline{u}%
_{i}(z_{0})e^{-qt_{0}}+\sum_{n=N_{t_{0}}}^{Nt}e^{-q\tau_{n}}p(l-Z_{\tau
_{n}^{-}}^{i}+Y_{n})1_{\{Z_{\tau_{n}^{-}}^{i}-Y_{n}-l<0\}}\\
&  & -\int_{t_{0}}^{t}e^{-qs}\left(  \sigma_{i}\overline{u}_{i}^{\prime}%
(Z_{s}^{i})-(q+\lambda)\overline{u}_{i}(Z_{s}^{i})+\lambda\int_{0}^{Z_{s^{-}%
}^{i}-l}\overline{u}_{i}(Z_{s^{-}}^{i}-\alpha)dF(\alpha)\right)  ds\\
&  & -\int_{t_{0}}^{t}e^{-qs}\left(  \lambda\int_{Z_{s^{-}}^{i}-l}^{\infty
}\left(  p(\alpha-Z_{s^{-}}^{i}+l)+\overline{u}_{i}(l)\right)  dF(\alpha
)\right)  ds
\end{array}
\label{Martingalas Mi}%
\end{equation}
with%
\[
\text{ }Z_{t}^{i}=z_{0}+\sigma_{i}\left(  t-t_{0}\right)  -\sum
\nolimits_{n=N_{t_{0}}}^{N_{t}}\min\{Y_{n},Z_{\tau_{n}^{-.}}^{i}-l\}\text{ for
}t\geq t_{0}\geq0,
\]
it can be seen that $M^{i}(z_{0},t_{0},t)$ is a martingale with zero
expectation for $t\geq t_{0}$.

Consider first the case $J_{k}=i$ and $J_{k+1}=j$ with $i=1,2$, $j=0,1,2$ and
$i\neq j$. Since $\overline{u}_{i}$ is absolutely continuous, the function
$t\rightarrow\overline{u}_{i}(X_{t}^{\pi})e^{-qt}$ is also absolutely
continuous, between the customer demands. Using an extension of the Dynkin's
Formula, we obtain
\[%
\begin{array}
[c]{l}%
\overline{u}_{j}(X_{T_{k+1}}^{\pi})e^{-qT_{k+1}}-\overline{u}_{i}(X_{T_{k}%
}^{\pi})e^{-qT_{k}}\\%
\begin{array}
[c]{ll}%
= & \overline{u}_{j}(X_{T_{k+1}}^{\pi})e^{-qT_{k+1}}-\overline{u}%
_{i}(X_{T_{k+1}}^{\pi})e^{-qT_{k+1}}+\overline{u}_{i}(X_{T_{k+1}}^{\pi
})e^{-qT_{k+1}}-\overline{u}_{i}(X_{T_{k}}^{\pi})e^{-qT_{k}}\\
\geq & -K_{ij}e^{-qT_{k+1}}+\overline{u}_{i}(X_{T_{k+1}}^{\pi})e^{-qT_{k+1}%
}-\overline{u}_{i}(X_{T_{k}}^{\pi})e^{-qT_{k}}\\
= & -K_{ij}e^{-qT_{k+1}}+\int_{T_{k}}^{T_{k+1}}e^{-qs}\mathcal{L}%
_{i}(\overline{u}_{i})(X_{s}^{\pi})ds\\
& -\left(  \int_{T_{k}}^{T_{k+1}}e^{-qs}h_{i}(X_{s}^{\pi})ds+\sum_{n=N_{T_{k}%
}}^{N_{T_{k+1}}}e^{-q\tau_{n}}p(l-X_{\tau_{n}^{-}}^{\pi}+Y_{n})1_{\{X_{\tau
_{n}^{-}}^{\pi}-Y_{n}-l<0\}}\right) \\
& +M^{i}(X_{T_{k}}^{\pi},T_{k},T_{k+1});
\end{array}
\end{array}
\]
and so, since $\overline{u}_{i}$ is a supersolution of (\ref{HJBi}), we get
that%
\[%
\begin{array}
[c]{l}%
\mathbb{E}\left[  \left.  \overline{u}_{j}(X_{T_{k+1}}^{\pi})e^{-qT_{k+1}%
}-\overline{u}_{i}(X_{T_{k}}^{\pi})e^{-qT_{k}}\right\vert \mathcal{F}_{T_{k}%
}\right] \\%
\begin{array}
[c]{cc}%
\geq & -\mathbb{E}\left[  K_{ij}e^{-qT_{k+1}}+\left.  \int_{T_{k}}^{T_{k+1}%
}e^{-qs}h_{i}(X_{s}^{\pi})ds+\sum_{n=N_{T_{k}}}^{N_{T_{k+1}}}e^{-q\tau_{n}%
}p(l-X_{\tau_{n}^{-}}^{\pi}+Y_{n})1_{\{X_{\tau_{n}^{-}}^{\pi}-Y_{n}%
-l<0\}}\right\vert \mathcal{F}_{T_{k}}\right]  .
\end{array}
\end{array}
\]

In the case $J_{k}=0$ we have $\mathcal{J}_{T_{k+1}}\neq0$, then $X_{s}^{\pi
}=b$ in $[T_{k},T_{k+1})$ , $T_{k+1}=\tau_{n}$ for some $n$ and so,
analogously to the previous case,
\[%
\begin{array}
[c]{l}%
\overline{u}_{\mathcal{J}_{T_{k+1}}}(X_{T_{k+1}}^{\pi})e^{-qT_{k+1}}%
-\overline{u}_{0}(X_{T_{k}}^{\pi})e^{-qT_{k}}\\%
\begin{array}
[c]{ll}%
= & e^{-qT_{k}}\left(  (\overline{u}(b-Y_{n})1_{\{b-Y_{n}-l\geq0\}}%
+\overline{u}_{\mathcal{J}_{T_{k+1}}}(l)1_{\{b-Y_{n}-l<0\}})e^{-q(T_{k+1}%
-T_{k})}-\overline{u}_{0}(b)\right)  -K_{0\mathcal{J}_{T_{k+1}}}e^{-qT_{k+1}%
}\\
= & e^{-qT_{k}}\left(  \overline{u}(b-Y_{n})1_{\{b-Y_{n}-l\geq0\}}%
+\overline{u}_{\mathcal{J}_{T_{k+1}}}(l)1_{\{b-Y_{n}-l<0\}}+p(l-b+Y_{n}%
)1_{\{b-Y_{n}-l<0\}}\right)  e^{-q(T_{k+1}-T_{k})}\\
& -e^{-qT_{k}}\left(  \frac{\lambda}{q+\lambda}\left(
{\textstyle\int_{0}^{b-l}}
\overline{u}(b-\alpha)dF(\alpha)+%
{\textstyle\int_{b-l}^{\infty}}
\left(  p(\alpha-b+l)+\overline{u}(l)\right)  dF(\alpha)\right)  \right) \\
& -\left(  K_{0\mathcal{J}_{T_{k+1}}}e^{-qT_{k+1}}+\int_{T_{k}}^{T_{k+1}%
}e^{-qs}h_{0}(b)ds+e^{-qT_{k+1}}p(l-b+Y_{n})1_{\{b-Y_{n}-l<0\}}\right) \\
& +\int_{T_{k}}^{T_{k+1}}e^{-qs}h_{0}(b)ds-\frac{\lambda}{q+\lambda}%
e^{-qT_{k}}h_{0}(b)\text{,}%
\end{array}
\end{array}
\]
and, since $T_{k+1}-T_{k}$ is distributed as $exp(\lambda)$, we obtain that%

\[%
\begin{array}
[c]{ll}%
0= & \mathbb{E}\left[  \left.  e^{-qT_{k}}\left(  \overline{u}(b-Y_{n}%
)1_{\{b-Y_{n}-l\geq0\}}+\overline{u}_{\mathcal{J}_{T_{k+1}}}(l)1_{\{b-Y_{n}%
-l<0\}}\right)  \right\vert \mathcal{F}_{T_{k}}\right] \\
& +\mathbb{E}\left[  \left.  e^{-qT_{k+1}}p(l-b+Y_{n})1_{\{b-Y_{n}%
-l<0\}}\right\vert \mathcal{F}_{T_{k}}\right] \\
& -\mathbb{E}\left[  \left.  \frac{e^{-qT_{k}}}{q+\lambda}\lambda\left(
{\textstyle\int_{0}^{b-l}}
\overline{u}(b-\alpha)dF(\alpha)+%
{\textstyle\int_{b-l}^{\infty}}
\left(  p(\alpha-b+l)+\overline{u}(l)\right)  dF(\alpha)\right)  \right\vert
\mathcal{F}_{T_{k}}\right]  .
\end{array}
\]
and so%
\[%
\begin{array}
[c]{l}%
\mathbb{E}\left[  \left.  \overline{u}_{\mathcal{J}_{T_{k+1}}}(X_{T_{k+1}%
}^{\pi})e^{-qT_{k+1}}-\overline{u}_{0}(X_{T_{k}}^{\pi})e^{-qT_{k}}\right\vert
\mathcal{F}_{T_{k}}\right] \\%
\begin{array}
[c]{cc}%
= & -\mathbb{E}\left[  \left.  K_{0\mathcal{J}_{T_{k+1}}}e^{-qT_{k+1}}%
+\int_{T_{k}}^{T_{k+1}}e^{-qs}h_{0}(b)ds+e^{-qT_{k+1}}p(l-b+Y_{n}%
)1_{\{b-Y_{n}-l<0\}}\right\vert \mathcal{F}_{T_{k}}\right]  .
\end{array}
\end{array}
\]
Analogously, we can prove that
\[%
\begin{array}
[c]{c}%
\begin{array}
[c]{l}%
\mathbb{E}\left[  \left.  \overline{u}_{J_{m_{t}}}(X_{t}^{\pi})e^{-qt}%
-\overline{u}_{_{J_{m_{t}}}}(X_{T_{m_{t}}}^{\pi})e^{-qT_{m_{t}}}\right\vert
\mathcal{F}_{T_{m_{t}}}\right] \\
\geq-\mathbb{E}\left[  \left.  \int_{T_{m_{t}}}^{t}e^{-qs}h_{_{J_{m_{t}}}%
}(X_{s}^{\pi})ds+\sum_{n=N_{T_{m_{t}}}}^{t}e^{-q\tau_{n}}p(l-X_{\tau_{n}^{-}%
}^{\pi}+Y_{n})1_{\{X_{\tau_{n}^{-}}^{\pi}-Y_{n}-l<0\}}\right\vert
\mathcal{F}_{T_{m_{t}}}\right]  .
\end{array}
\end{array}
\]
Taking the expected value in (\ref{Expresion Inicial}), we obtain
\[%
\begin{array}
[c]{lll}%
\mathbb{E}\left[  \overline{u}_{\mathcal{J}_{t}}(X_{t}^{\pi})e^{-qt}\right]
-\overline{u}_{j}(x) & = & \mathbb{E}\left[  \sum_{k=0}^{m_{t}-1}%
\mathbb{E}\left[  \left.  \left(  \overline{u}_{_{J_{k+1}}}(X_{T_{k+1}}^{\pi
})e^{-qT_{k+1}}-\overline{u}_{_{J_{k}}}(X_{T_{k}}^{\pi})e^{-qT_{k}}\right)
\right\vert \mathcal{F}_{T_{k}}\right]  \right] \\
&  & +\mathbb{E}\left[  \mathbb{E}\left[  \left.  (\overline{u}_{J_{m_{t}}%
}(X_{t}^{\pi})e^{-qt}-\overline{u}_{_{J_{m_{t}}}}(X_{T_{m_{t}}}^{\pi
})e^{-qT_{m_{t}}})\right\vert \mathcal{F}_{T_{m_{t}}}\right]  \right] \\
& \geq & -V_{j}^{\pi}(x)
\end{array}
\]
taking the limit with $t$ going to infinity, and using that $X_{t}^{\pi}%
\in\lbrack l,b]$ we obtain that $\overline{u}_{j}(x)\leq V_{j}^{\pi}(x)$ for
$j=1,2$.

Considering instead the controlled risk process $X_{t}^{\pi}$ starting at $b$,
we obtain with a similar proof that $\overline{u}_{0}(b)\leq V_{0}^{\pi}(b)$.
$\blacksquare$

From Propositions \ref{V is viscosity solution} and
\ref{Largest Viscosity Supersolution}, we obtain the following verification result.

\begin{theorem}
\label{Teorema Verificacion} Consider two families of admissible strategies
$\left\{  \pi_{x,i}\in\Pi_{x,i}:x\in\lbrack l,b)\right\}  $ for $i=1,2$. If
the functions $w_{i}(x):=V_{i}^{\pi_{x,i}}(x)$ for $i=1,2$ are viscosity
supersolutions of the respective HJB equation (\ref{HJBi}) for $x\in(l,b)$ and
satisfy the boundary conditions%
\[
w_{1}(b)\leq w_{0}(b)+K_{10},w_{2}(b)\leq w_{0}(b)+K_{20},
\]
where%
\[%
\begin{array}
[c]{l}%
w_{0}(b)=\frac{\lambda}{q+\lambda}\left(
{\textstyle\int_{0}^{b-l}}
\overline{w}(b-\alpha)dF(\alpha)+%
{\textstyle\int_{b-l}^{\infty}}
\left(  p(\alpha-b+l)+\overline{w}(l)\right)  dF(\alpha)\right)  +\frac
{h_{0}(b)}{q+\lambda}\text{ and}\\
\overline{w}(x)=\min\{K_{01}+w_{1}(x),K_{02}+w_{2}(x)\}.
\end{array}
\]
Then, $w_{0}(b)=V_{0}(b)$ and $w_{i}=V_{i}$ for $i=1,2.$

\end{theorem}

In the remainder of the section, we show that there exists an optimal
production-inventory strategy and it is stationary in the sense that depends
only on the phase and the inventory level.

\begin{definition}
\label{Estrategia Estacionaria}Given two disjoint closed sets $A_{12}$ and
$A_{21}$ in $[l,b)$ and a closed set $C_{1}$ in $[l,b)$ with $A_{21}\subset
C_{1}$ and $A_{12}\subset\lbrack l,b)-C_{1}$, we define the
production-inventory band strategy associated to the sets $(A_{12}%
,A_{21},C_{1})$ as follows:
\end{definition}

\begin{enumerate}
\item If the current phase is $i=1$ and the current inventory level is $x\in
A_{12}$, change immediately to phase $2$, if the current inventory level
$x\in\lbrack l,b)-A_{12}$ stay in phase 1.

\item If the current phase is $i=2$ and the current inventory level is $x\in
A_{21}$, change immediately to phase $1$, if the current inventory level
$x\in\lbrack l,b)-A_{21}$ stay in phase 2.

\item If the current phase is $i=0$ with current inventory level $b$, then in
the event of an arrival of the next customer demand of size $Y,$ switch on the
production to phase $1$ if $\max\{b-Y,l\}\in C_{1}$ and switch on the
production to phase $2$ if $\max\{b-Y,l\}\in\lbrack l,b)-C_{1}$.

\item If the inventory level reaches $b$, it is mandatory to switch to phase
$0$.
\end{enumerate}

The sets $A_{ij}$ are called the \textit{switching zone} from the phase $i$ to
phase $j$, and the sets $C_{1}$ and $C_{2}=[l,b)-C_{1}$ are called the
\textit{selection zones} for phases 1 and 2 respectively. Also, the set
$[l,b)-\left(  A_{12}\cup A_{21}\right)  $ is called the \textit{non-action
zone}.

\begin{remark}
\label{Punto Fijo} Given the sets $\mathcal{A}=(A_{12},A_{21},C_{1})$, an
initial inventory level $x$ and an initial phase $i,$ we define and admissible
strategy $\pi_{x,i}^{\mathcal{A}}=(T_{k},J_{k})_{k\geq0}\in\Pi_{x,i}$ where
$J_{0}=i$ and $T_{k}$ is the $k$-th switching (from regime $J_{k-1}$ to
$J_{k})$ given by $(1)$, $(2)$, $(3)$ and $(4)$. Note that the switching times
$T_{k}$ are the times in which the controlled inventory process in $[l,b]$
exit the sets $[l,b)-A_{12},~[l,b)-A_{21}$ and $\{b\}$. Let us denote the cost
function of this admissible strategies as
\[
W_{i}^{\mathcal{A}}(x)=V_{i}^{\pi_{x,i}^{\mathcal{A}}}(x)\text{ for
}i=1,2\text{ and }x\in\lbrack l,b)\text{; and }W_{0}^{\mathcal{A}}%
(b)=V_{i}^{\pi_{b,0}^{\mathcal{A}}}(b).
\]
We can characterize the triple $(W_{0}^{\mathcal{A}}(b),W_{1}^{\mathcal{A}%
},W_{2}^{\mathcal{A}})$ as the unique fixed point of a contraction operator:
Let $\mathcal{C}[l,b)$ be the set of all the functions $W:[l,b)\rightarrow
\mathbf{R}$ continuous and bounded and let consider the Banach space
\[
\mathcal{B}=\mathbf{R}\times\mathcal{C}[l,b)\times\mathcal{C}[l,b)
\]
with norm
\[
\left\Vert (f_{0},f_{1},f_{2})\right\Vert =\max\{\left\vert f_{0}\right\vert
,\sup_{x\in\lbrack l,b)}\left\vert f_{1}(x)\right\vert ,\sup_{x\in\lbrack
l,b)}\left\vert f_{2}(x)\right\vert \}.
\]
We define, the operator $\mathcal{T}^{\mathcal{A}}:\mathcal{B\rightarrow B}$
as
\begin{equation}
\mathcal{T}^{\mathcal{A}}(f_{0},f_{1},f_{2})=\left(  \mathcal{T}%
_{0}^{\mathcal{A}}(f_{1},f_{2},f_{0}),\mathcal{T}_{1}^{\mathcal{A}}%
(f_{1},f_{2},f_{0}),\mathcal{T}_{2}^{\mathcal{A}}(f_{1},f_{2},f_{0})\right)  .
\label{TA}%
\end{equation}
We define $\mathcal{T}_{0}^{\mathcal{A}}$ as
\[%
\begin{array}
[c]{lll}%
\mathcal{T}_{0}^{\mathcal{A}}(f_{0},f_{1},f_{2}) & := & \mathbb{E}\left[
\int_{0}^{\tau_{1}}e^{-qs}h_{0}(b)ds)\right] \\
&  & +\mathbb{E}\left[  1_{\left\{  Y_{1}\leq b-l\right\}  }e^{-q\tau_{1}%
}\left(  \overline{f}(b-Y_{1})\right)  \right] \\
&  & +\mathbb{E}\left[  1_{\left\{  Y_{1}>b-l\right\}  }e^{-q\tau_{1}}\left(
p(Y_{1}-b+l)+\overline{f}(l)\right)  \right]
\end{array}
\]
where
\[%
\begin{array}
[c]{lll}%
\overline{f}(x) & := & (f_{1}(x)+K_{01})1_{\left\{  x\in C_{1}^{\ast}\right\}
}+(f_{2}(x)+K_{02})1_{\left\{  x\notin C_{1}^{\ast}\right\}  },
\end{array}
\]
here $(\tau_{1},Y_{1})$ is the time and size of the first costumer demand.
Take the admissible strategy $\pi_{x,i}^{\mathcal{A}}=(T_{k},J_{k})_{k\geq
0}\in\Pi_{x,i}$ as defined in Definition \ref{Estrategia Estacionaria} and
consider the associated controlled inventory process $X_{t}$ and the process
$\mathcal{J}_{t}\ $defined in (\ref{Jt}), we define $\mathcal{T}%
_{i}^{\mathcal{A}}$ as%
\[%
\begin{array}
[c]{lll}%
\mathcal{T}_{i}^{\mathcal{A}}(f_{0},f_{1},f_{2})(x) & = & \mathbb{E}\left[
\int_{0}^{\tau_{1}}e^{-qs}h_{\mathcal{J}s}(X_{s})ds\right]  +\\
&  & +\mathbb{E}\left[
{\textstyle\sum\nolimits_{k=1}^{\infty}}
1_{\left\{  T_{k}<\tau_{1}\right\}  }e^{-qT_{k}}K_{J_{k-1}J_{k}}\right] \\
&  & +\mathbb{E}\left[  1_{\left\{  X_{\tau_{1}^{-}}-Y_{1}\geq l\right\}
}e^{-q\tau_{1}}\left(  f_{\mathcal{J}\tau_{1}^{-}}(X_{\tau_{1}^{-}}%
-Y_{1})\right)  \right] \\
&  & +\mathbb{E}\left[  1_{\left\{  X_{\tau_{1}^{-}}-Y_{1}<l\right\}
}e^{-q\tau_{1}}\left(  p(l-X_{\tau_{1}^{-}}+Y_{1})+f_{\mathcal{J}\tau_{1}^{-}%
}(l)\right)  \right]
\end{array}
\]
for $x\in\lbrack l,b)$ and $i\in\{1,2\}.$ Note that
\[%
\begin{array}
[c]{ccc}%
\left\vert \mathcal{T}_{i}^{\mathcal{A}}(f_{0},f_{1},f_{2})-\mathcal{T}%
_{i}^{\mathcal{A}}(g_{0},g_{1},g_{2})\right\vert  & \leq & \mathbb{E}\left(
e^{-q\tau_{1}}\right)  \left\Vert (f_{0},f_{1},f_{2})-(g_{0},g_{1}%
,g_{2})\right\Vert \\
& = & \frac{\lambda}{q+\lambda}\left\Vert (f_{0},f_{1},f_{2})-(g_{0},g_{1}%
,g_{2})\right\Vert
\end{array}
\]
and so $\mathcal{T}^{\mathcal{A}}:\mathcal{B\rightarrow B}$ is a contraction
operator with a unique fixed point. Finally, by the definition of the
production-inventory strategy associated to the sets $\mathcal{A=}%
(A_{12},A_{21},C_{1})$, it follows immediately that the triple $(W_{0}%
^{\mathcal{A}}(b),W_{1}^{\mathcal{A}},W_{2}^{\mathcal{A}})$ is a fixed point
of the operator $\mathcal{T}^{\mathcal{A}}$.
\end{remark}

In the following theorem we prove that there exists an optimal strategy and
that it comes from a production-inventory band strategy as defined in
Definition \ref{Estrategia Estacionaria}.

\begin{theorem}
The optimal strategy of problem (\ref{Definicion Vi}) and (\ref{Definicion V0}%
), is the production-inventory strategy associated to the sets $\mathcal{A}%
^{\ast}=(A_{12}^{\ast},A_{21}^{\ast},C_{1}^{\ast})$ where%
\[%
\begin{array}
[c]{lll}%
A_{12}^{\ast} & = & \{x\in\lbrack l,b):V_{2}(x)+K_{12}-V_{1}(x)=0\},\\
A_{21}^{\ast} & = & \{x\in\lbrack l,b):V_{1}(x)+K_{21}-V_{2}(x)=0\},\\
C_{1}^{\ast} & = & \{x\in\lbrack l,b):K_{01}+V_{1}(x)\leq K_{02}+V_{2}(x)\}.
\end{array}
\]

\end{theorem}

Proof.

By Remark \ref{Punto Fijo}, it is enough to prove that the triple $\left(
V_{0}(b),V_{1},V_{2}\right)  $ is a fixed point of the operator $\mathcal{T}%
^{\mathcal{A}^{\ast}}\ $for the sets $\mathcal{A}^{\ast}=(A_{12}^{\ast}%
,A_{21}^{\ast},C_{1}^{\ast})$. By definition of the sets $A_{ij}^{\ast}$ and
$C_{1}^{\ast},$ we obtain immediately that $\mathcal{T}_{0}^{\mathcal{A}%
^{\ast}}(V_{0}(b),V_{1},V_{2})=V_{0}(b).$ Let us prove now that $\mathcal{T}%
_{i}^{\mathcal{A}^{\ast}}(V_{0}(b),V_{1},V_{2})(x)=V_{i}(x)$ for $x\in\lbrack
l,b)$ and $i=1,2.$ Since $\mathcal{L}_{0}(V_{0})(b)=0;$ and for
$\{i,j\}=\{1,2\}$ the functions $t\rightarrow V_{i}(X_{t})$ are absolutely
continuos, $\mathcal{L}_{i}(V_{i})=0$ a.e. in $[l,b)-A_{ij}^{\ast}$ and
$V_{j}(x)+K_{ij}-V_{i}(x)=0$ in $A_{ij}^{\ast}$; we can prove, with arguments
similar to the proof of Proposition \ref{Largest Viscosity Supersolution}, and
using the martingales introduced in (\ref{Martingalas Mi}) that%

\[%
\begin{array}
[c]{lll}%
\mathcal{T}_{i}^{\mathcal{A}^{\ast}}(V_{0},V_{1},V_{2})(x)-V_{i}(x) & = &
\mathbb{E}\left[  \int_{0}^{\tau_{1}}e^{-qs}h_{\mathcal{J}s}(X_{s})ds\right]
+\\
&  & +\mathbb{E}\left[
{\textstyle\sum\nolimits_{k=1}^{\infty}}
1_{\left\{  T_{k}<\tau_{1}\right\}  }e^{-qT_{k}}K_{J_{k-1}J_{k}}\right] \\
&  & +\mathbb{E}\left[  1_{\left\{  X_{\tau_{1}^{-}}-Y_{1}\geq l\right\}
}e^{-q\tau_{1}}\left(  V_{\mathcal{J}\tau_{1}^{-}}(X_{\tau_{1}^{-}}%
-Y_{1})\right)  \right] \\
&  & +\mathbb{E}\left[  1_{\left\{  X_{\tau_{1}^{-}}-Y_{1}<l\right\}
}e^{-q\tau_{1}}\left(  p(l-X_{\tau_{1}^{-}}+Y_{1})+V_{\mathcal{J}\tau_{1}^{-}%
}(l)\right)  \right]  -V_{i}(x).\\
& = & \mathbb{E}\left[  \int_{0}^{\tau_{1}}e^{-qs}\mathcal{L}_{\mathcal{J}%
s}(V_{\mathcal{J}s})(X_{s})ds\right] \\
& = & 0.
\end{array}
\]

Hence, $W_{0}^{\mathcal{A}^{\ast}}(b)=V_{0}(b),W_{1}^{\mathcal{A}^{\ast}%
}=V_{1}$ and $W_{2}^{\mathcal{A}^{\ast}}=V_{2}$.$\blacksquare$

\section{Finite Band strategies \label{Section Finite Band Strategies}}

We define the finite band strategies as the production-inventory band
strategies in which the non-action set $[l,b)-\left(  A_{{\small 12}}\cup
A_{{\small 21}}\right)  $ has a finite number of connected components.

Doshi et al. \cite{Doshi} studied the production-inventory band strategies
with switching zones $A_{{\small 12}}=[y_{{\small 1}},b)$\ and $A_{{\small 21}%
}=[l,y_{{\small 2}}]$ and selection zones $C_{{\small 1}}=[l,y_{{\small 2}}]$
and $C_{{\small 2}}=(y_{{\small 2}},b)$ for $l\leq y_{{\small 2}%
}<y_{{\small 1}}<b$.

Assuming that the optimal strategy is a finite band strategy, we look for it
in the following way;

\textbf{First step}. We find the best \textit{Doshi strategy}, that is we
construct the cost functions $(W_{0}^{\mathcal{A}}(b),W_{1}^{\mathcal{A}%
},W_{2}^{\mathcal{A}})$ for $\mathcal{A=}\left(  [y_{{\small 1}}%
,b),[l,y_{{\small 2}}],(y_{{\small 2}},b)\right)  $; then we minimize the
$W_{0}^{\mathcal{A}}(b)$ among the two variables $l\leq y_{{\small 2}%
}<y_{{\small 1}}<b$. We check whether the associated cost functions
$W_{0}^{\mathcal{A}}(b),W_{1}^{\mathcal{A}}$ and $W_{2}^{\mathcal{A}}$ of this
strategy satisfy the conditions of Theorem \ref{Teorema Verificacion}, if they
do this is the optimal strategy; if this is not the case, we go to the second step.

\textbf{Second step}. We consider the \textit{band strategies of type one}
where the non-action zone has one connected component. Here, the switching
zones are of the form $A_{{\small 12}}=[y_{{\small 1}},b)$\ and
$A_{{\small 21}}=[l,y_{{\small 2}}]$ and the selection zones are of the form
$C_{{\small 1}}=[l,y_{{\small 3}}]$ and $C_{{\small 2}}=(y_{{\small 3}},b)$
for $l\leq y_{{\small 2}}\leq y_{3}<y_{{\small 1}}<b$; the non-action zone is
$(y_{2},y_{1})$. Then we minimize $W_{0}^{\mathcal{A}}(b)$ among the three
variables $y_{{\small 2}},~y_{3}\,,~y_{{\small 1}}$. As before, we check
whether the associated cost functions $W_{0}^{\mathcal{A}}(b),W_{1}%
^{\mathcal{A}}$ and $W_{2}^{\mathcal{A}}$ of this strategy satisfy the
conditions of Theorem \ref{Teorema Verificacion}, if they do this is the
optimal strategy; if this is not the case, we go to the third step. Note that
the Doshi strategies are the band strategies of type one in which
$y_{{\small 2}}=y_{3}.$

\textbf{Third step}. We consider the \textit{band strategies of type two}
where the non-action zone has two connected components. Here, the switching
zones are of the form $A_{{\small 12}}=[y_{{\small 1}},y_{4}]$\ and
$A_{{\small 21}}=[l,y_{{\small 2}}]$ and the selection zones are of the form
$C_{{\small 1}}=[l,y_{{\small 3}}]$ and $C_{{\small 2}}=(y_{{\small 3}},b)$
for $l\leq y_{{\small 2}}\leq y_{3}<y_{{\small 1}}<y_{{\small 4}}<b$; in these
band strategies, the non-action zone $(y_{{\small 2}},y_{{\small 1}}%
)\cup(y_{{\small 4}},b)$ has two connected components. Now, we minimize
$W_{0}^{\mathcal{A}}(b)$ among the four variables $y_{{\small 2}}%
,~y_{3},~y_{{\small 1}},~y_{{\small 4}}$. Again, we check whether the
associated cost functions $W_{0}^{\mathcal{A}}(b),W_{1}^{\mathcal{A}}$ and
$W_{2}^{\mathcal{A}}$ of this strategy satisfy the conditions of Theorem
\ref{Teorema Verificacion}, if they do this is the optimal strategy; if this
is not the case, we consider band strategies where the non-action zone has
more connected components. And so on...

In the next section, we describe how to find the cost functions of band
strategies with one and two connected components using scale functions. We
also show how to find the decomposition into the different types of costs:
holding, production, switching and penalty costs.

In Section \ref{Section Examples}, we show examples where the optimal strategy
is a Doshi strategy (Figure 1), is a band strategy of type one (Figure 6) and
is a band strategy of type two (Figure 11).

\section{The value functions of band strategies
\label{Section Value Functions}}

In this section we derive the cost functions for band strategies of type one and two. Throughout this section we assume that $l=0$.
We further assume that the holding cost per time unit in phase $i$ when the
inventory level is $x$ is : $h_{i}(x)=a_{i}+c_{i}x\ $for $\,i=1,2$, where
$a_{i},c_{i}\geq0$ are given. To obtain the value function we apply the
fluctuation theory for L\'{e}vy processes as described in Chapter 8 in Kyprianou
(2014) and Avram et al. (2019).

\subsection{Preliminaries}

For $i=1,2$ let
\[
X_{i,t}:=x+\sigma_{i}t-\sum_{n=1}^{N(t)}Y_{n}%
\]
be the uncontrolled process at phase $i$ with initial inventory level $x$. The
processes $X_{i}$ are spectrally negative bounded variation L\'{e}vy processes.
Let us define
\[
\varphi_{i}(\theta)=\log\mathbb{E}\left[  e^{\theta\left(  X_{i,1}-x\right)
}\right]  =\sigma_{i}\theta-\lambda+\lambda\mathcal{L}_{Y}(\theta),
\]
where $\mathcal{L}_{Y}(\theta):=\mathbb{E}[e^{-\theta Y_{1}}]$. Let us also
define the exit times $\tau_{i,a}^{-}=\inf\{t:X_{i,t}<a\}$ and $\tau_{i,d}%
^{+}=\inf\{t:X_{i,t}=d\}$.

In this section, we use the following notations:

\begin{description}
\item[$\bullet$] $W_{i}^{(q)}(x)$ (the scale function associated with $X_{i}%
$). This scale function is defined by its Laplace transform:
\[
\int_{0}^{\infty}e^{-\theta x}W_{i}^{(q)}(x)dx=\frac{1}{q-\varphi_{i}(\theta
)}.
\]

\item[$\bullet$]
\[
Z_{i}^{(q)}(x,\theta)=e^{\theta x}\left(  1+(q-\varphi_{i}(\theta))\int
_{0}^{x}e^{-\theta y}W_{i}^{(q)}(y)dy\right)  .
\]

Denote $Z_{i}^{(q)}(x)=Z_{i}^{(q)}(x,0)=1+q\int_{0}^{x}W^{(q)}_{i}(y)dy$.

\item[$\bullet$] $\overline{W}_{i}^{(q)}(x)=\int_{0}^{x}W_{i}^{(q)}(y)dy.$

\item[$\bullet$] $\overline{\overline{W}}_{i}^{(q)}(x)=\int_{0}^{x}%
\overline{W}_{i}^{(q)}(y)dy.$

\item[$\bullet$] $\overline{Z}_{i}^{(q)}(x)=\int_{0}^{x}Z_{i}^{(q)}%
(y)dy=x+q\overline{\overline{W}}_{i}(x)$.
\end{description}

Throughout this section, we will also use the following results:

\begin{description}
\item[$\bullet$]
\begin{equation}
\mathbb{E}_{x}\left[  e^{-q\tau_{i,d}^{+}}1_{\tau_{i,d}^{+}<\tau_{i,a}^{-}%
}\right]  =\frac{W_{i}^{(q)}(x-a)}{W_{i}^{(q)}(d-a)}. \label{eq.exitup}%
\end{equation}

\item[$\bullet$]
\begin{align}
&  \mathbb{E}_{x}\left[  e^{-q\tau_{i,0}^{-}+\theta X_{i,\tau_{i,0}^{-}}%
}1_{\tau_{i,0}^{-}<\tau_{i,d}^{+}}\right] \nonumber\label{eq.exitdeficit}\\
&  =Z_{i}^{(q)}(x,\theta)-\frac{W_{i}^{(q)}(x)}{W_{i}^{(q)}(d)}Z_{i}%
^{(q)}(d,\theta).
\end{align}

\item[$\bullet$]
\begin{align}
&  \mathbb{E}_{x}\left[  e^{-q\tau_{i,0}^{-}}1_{\tau_{i,0}^{-}<\tau_{i,d}^{+}%
}\right] \nonumber\label{eq.exitdown}\\
&  =Z_{i}^{(q)}(x)-\frac{W_{i}^{(q)}(x)}{W_{i}^{(q)}(d)}Z_{i}^{(q)}(d).
\end{align}

\item[$\bullet$] For $0<y<d,$ let us define the $q$-potential measure of
$X_{i}$ as
\begin{equation}
U_{i}^{(q)}(a,d,x,dy)=\int_{0}^{\infty}e^{-qt}\mathbb{P}_{x}\left[  X_{i,t}\in
dy,\tau_{i,a}^{-}\wedge\tau_{i,d}^{+}>t\right]  dt. \label{eq.U}%
\end{equation}

\item By Theorem 8.7 in Kyprianou (2014), $U_{i}^{(q)}(a,d,x,dy)=u_{i}%
^{(q)}(x,y)dy,$ where
\begin{equation}
u_{i}^{(q)}(a,d,x,y)=\frac{W_{i}^{(q)}(x-a)W_{i}^{(q)}(d-y)}{W_{i}^{(q)}%
(d-a)}-W_{i}^{(q)}(x-y). \label{eq.u}%
\end{equation}

\end{description}

Throughout, we denote by $\mathbb{E}^{i}$ the expectation according to the
probability law $\mathbb{P}^{i}$ induced by the process $X_{i}$ for $i=1,2$.

\subsection{Cost functions for strategies of type one.}

As defined in the previous sections the switching zone are $A_{{\small 12}%
}=[y_{{\small 1}},b)$\ and $A_{{\small 21}}=[0,y_{{\small 2}}]$ and the
selection zones are $C_{{\small 1}}=[0,y_{{\small 3}}]$ and $C_{{\small 2}%
}=(y_{{\small 3}},b)$ for $0\leq y_{{\small 2}}\leq y_{3}<y_{{\small 1}}<b$;
the non-action zone is $(y_{2},y_{1})$. The value function is obtained in
three steps, first we obtain the expected discounted holding cost, then the
expected discounted shortage cost and finally the expected discounted
switching cost.

\subsubsection{Expected discounted holding cost.}

Here, we compute the formulas for:

\begin{itemize}
\item[$\bullet$] $\mathcal{H}_{i}^{(q)}(x)$-- the expected discounted holding
cost starting at $x$ at phase $i$, $i=1,2$.

\item[$\bullet$] $\mathcal{H}_{0}^{(q)}(b)$-- the expected discounted holding
cost starting at $b$.

\item[$\bullet$] $H_{1}^{(q)}(x,y_{1})$-- the expected discounted holding cost
until reaching $y_{1}$ starting at $x$ at phase $1$, $0\leq x<y_{1}$.

\item[$\bullet$] $H_{2}^{(q)}(x,y_{2},b)$-- the expected discounted holding
cost until reaching $b$ or down-crossing $y_{2}$ starting at $x$ at phase $2$,
$y_{2}<x<b$.
\end{itemize}

In order to do that, let us define $\underline{X}_{1,t}=\inf\{s\leq
t,X_{1,s}\}$ and $L_{t}=-(\underline{X}_{1,t}\wedge0)$. Let $R_{t}%
=X_{1,t}+L_{t}$. Let $\kappa_{y_{1}}^{+}=\inf\left\{  {t:R_{t}\geq y_{1}%
}\right\}  $ be the first time that $R$ reaches $y_{1}$. Notice that when the
inventory is less than $y_{1}$ and the phase is $1$, the inventory evolves as
$R$. By Theorem 8.1 (ii) in Kyprianou (2014),
\begin{equation}
\mathbb{E}^{1}_{x}[e^{-q\kappa_{y_{1}}^{+}}]=\frac{Z_{1}^{(q)}(x)}{Z_{1}%
^{(q)}(y_{1})}. \label{kypkappa}%
\end{equation}

\begin{remark}
The main tool to evaluate the expected discounted holding cost is the
Kella--Whitt martingale, \cite{kella}: let $X_{t}$ be a spectrally negative L\'{e}vy process
with Laplace exponent $\varphi(\alpha)=\log\mathbb{E}[e^{\alpha\left(
X_{1}-X_{0}\right)  }]$, $Y_{t}$ an adapted process with bounded expected
variation on finite intervals and $V_{t}=X_{t}+Y_{t}$. Let $\Delta Y_{s}%
=Y_{s}-Y_{s^{-}}$ and $Y^{c}$ the continuous part of $Y$, i.e. $Y_{t}%
^{c}=Y_{t}-\sum_{0\leq s\leq t}\Delta Y_{s}$. Then:
\begin{align}
&  M_{t}=\varphi(\alpha)\int_{0}^{t}e^{\alpha V_{s}}ds+e^{\alpha V_{0}%
}-e^{\alpha V_{t}}+\alpha\int_{0}^{t}e^{\alpha V_{s}}dY_{s}^{c}%
\nonumber\label{eq.kella-whitt}\\
&  +\sum_{0\leq s\leq t}e^{\alpha V_{s}}(1-e^{-\alpha\Delta Y_{s}})
\end{align}
is a zero mean martingale.
\end{remark}

From the strong Markov property at $y_{1}$ and (\ref{kypkappa}), it follows
that
for $0<x<y_{1}$,
\begin{equation}
\mathcal{H}_{1}^{(q)}(x)=H_{1}^{(q)}(x,y_{1})+\mathbb{E}_{x}^{1}%
[e^{-q\kappa_{y_{1}}^{+}}]\mathcal{H}_{2}^{(q)}(y_{1})=H_{1}^{(q)}%
(x,y_{1})+\frac{Z_{1}^{(q)}(x)}{Z_{1}^{(q)}(y_{1})}\mathcal{H}_{2}^{(q)}%
(y_{1}). \label{eq.calH1}%
\end{equation}

Similarly, for $y_{2}<x<b$,
\begin{align}
&  \mathcal{H}_{2}^{(q)}(x)=H_{2}^{(q)}(x,y_{2},b)+\mathbb{E}_{x}^{2}\left[
e^{-q\tau_{2,y_{2}}^{-}}1_{\tau_{2,y_{2}}^{-}<\tau_{2,b}^{+}}\mathcal{H}%
_{1}^{(q)}(X_{2,\tau_{2,y_{2}}^{-}})\right] \nonumber\\
&  +\mathbb{E}_{x}\left[  e^{-q\tau_{b}^{+}}1_{\tau_{2,b}^{+}<\tau_{2,y_{2}%
}^{-}}\right]  \mathcal{H}_{0}^{(q)}(b). \label{eq.calH2}%
\end{align}
And, if we denote by $Exp(\lambda)$ an exponentially distributed random
variable with rate $\lambda$,%

\begin{align}
&  \mathcal{H}_{0}^{(q)}(b)=h_{0}(b)\mathbb{E}\left[  \int_{0}^{Exp(\lambda
)}e^{-qt}dt\right] \nonumber\label{eq.calHb}\\
&  +\mathbb{E}\left[  e^{-qExp(\lambda)}\left(  \int_{0}^{b-y_{3}}%
\mathcal{H}_{2}^{(q)}(b-z)dF(z)+\int_{b-y_{3}}^{\infty}\mathcal{H}_{1}%
^{(q)}(b-z)dF(z)\right)  \right] \nonumber\\
&  =\frac{h_{0}(b)}{q+\lambda}+\frac{\lambda}{\lambda+q}\left(  \int
_{0}^{b-y_{3}}\mathcal{H}_{2}^{(q)}(b-z)dF(z)+\int_{b-y_{3}}^{\infty
}\mathcal{H}_{1}^{(q)}(b-z)dF(z)\right)  .
\end{align}

\begin{remark}
For $x<0,$ we define $H_{1}^{(q)}(x,y_{1})=H_{1}^{(q)}(0,y_{1})$ and
$\mathcal{H}_{1}^{(q)}(x)=\mathcal{H}_{1}^{(q)}(0)$.
\end{remark}

First, we find the formula for $H_{2}^{(q)}(x,y_{2},b)$:
\begin{align}
&  H_{2}^{(q)}(x,y_{2},b)=\mathbb{E}_{x}^{2}\left[  \int_{0}^{\tau_{2,b}%
^{+}\wedge\tau_{2,y_{2}}^{-}}e^{-qt}(a_{2}+c_{2}X_{2,t})dt\right]
\nonumber\label{eq.H2}\\
&  =\frac{a_{2}}{q}(1-\mathbb{E}_{x}^{2}[e^{-q(\tau_{2,b}^{+}\wedge
\tau_{2,y_{2}}^{-})}])+c_{2}\mathbb{E}_{x}^{2}[\int_{0}^{\tau_{2,b}^{+}%
\wedge\tau_{2,y_{2}}^{-}}e^{-qt}X_{2,t}dt].
\end{align}
Let
\begin{equation}
h_{2,1}(x,y_{2},b)=\frac{1}{q}\left(  1-\mathbb{E}_{x}^{2}[e^{-q(\tau
_{2,b}^{+}\wedge\tau_{2,y_{2}}^{-})}]\right)  .
\end{equation}
Applying (\ref{eq.exitdown}) and (\ref{eq.exitup}) yields:
\begin{align}
&  h_{2,1}(x,y_{2},b)=\frac{1}{q}\left(  1-\mathbb{E}_{x}^{2}\left[
e^{-q\tau_{2,y_{2}}^{-}}\,1_{\tau_{2,y_{2}}^{-}<\tau_{2,b}^{+}}\right]
-\mathbb{E}_{x}^{2}\left[  e^{-q\tau_{2,b}^{+}}\,1_{\tau_{2,b}^{+}%
<\tau_{2,y_{2}}^{-}}\right]  \right) \nonumber\\
&  =\frac{1}{q}\left(  1-Z_{2}^{(q)}(x-y_{2})+\frac{W_{2}^{(q)}(x-y_{2}%
)}{W_{2}^{(q)}(b-y_{2})}Z_{2}^{(q)}(b-y_{2})-\frac{W_{2}^{(q)}(x-y_{2})}%
{W_{2}^{(q)}(b-y_{2})}\right)  . \label{eq.h21}%
\end{align}



In order to obtain $\mathbb{E}_{x}^{2}[\int_{0}^{\tau_{2,b}^{+}\wedge
\tau_{2,y_{2}}^{-}}e^{-qt}X_{2,t}dt]$, we apply Kella-Whitt martingale
(\ref{eq.kella-whitt}) for $X_{2,t},$ $\varphi_{2}(\alpha)=\log\mathbb{E}%
[e^{\alpha\left(  X_{2,1}-X_{2,0}\right)  }]$ and $Y_{t}=-qt/\alpha$, so%

\begin{equation}
\mathbb{E}_{x}^{2}\left[  (\varphi_{2}(\alpha)-q)\int_{0}^{\tau_{2,b}%
^{+}\wedge\tau_{2,y_{2}}^{-}}e^{\alpha X_{2,s}-qs}ds+e^{\alpha x}-e^{\alpha
X_{2,\tau_{2,b}^{+}\wedge\tau_{2,y_{2}}^{-}}}\right]  =0. \label{eq.kellaH2}%
\end{equation}
Taking derivative of (\ref{eq.kellaH2}) with respect to $\alpha$ at
$\alpha=0,$ we obtain%

\begin{equation}
\varphi_{2}^{\prime}(0)\mathbb{E}_{x}^{2}\left[  \int_{0}^{\tau_{2,y_{2}}%
^{-}\wedge\tau_{2,b}^{+}}e^{-qs}ds\right]  +x-\frac{\partial}{\partial\alpha
}\mathbb{E}_{x}^{2}\left[  e^{\alpha X_{2,(\tau_{2,y_{2}}^{-}\wedge\tau
_{2,b}^{+})}-q(\tau_{2,y_{2}}^{-}\wedge\tau_{2,b}^{+})}|_{\alpha=0}\right]
=q\mathbb{E}_{x}^{2}\left[  \int_{0}^{\tau_{2,y_{2}}^{-}\wedge\tau_{2,b}^{+}%
}X_{2,s}e^{-qs}ds\right]  . \label{eq.jointdef}%
\end{equation}
By (\ref{eq.exitdeficit}), (\ref{eq.kellaH2}) and (\ref{eq.jointdef}), we get
\begin{align}
&  \mathbb{E}_{x}^{2}\left[  e^{\alpha X_{2,(\tau_{2,y_{2}}^{-}\wedge
\tau_{2,b}^{+})}-q(\tau_{2,y_{2}}^{-}\wedge\tau_{2,b}^{+})}\right]
\nonumber\label{eq.jointdef1}\\
&  =\mathbb{E}_{x}^{2}\left[  e^{\alpha b-q\tau_{2,b}^{+}}\,1_{\tau_{2,b}%
^{+}<\tau_{2,y_{2}}^{-}}\right]  +\mathbb{E}_{x}^{2}\left[  e^{\alpha
X_{2,\tau_{2,y_{2}}^{-}}-q\tau_{2,y_{2}}^{-}}\,1_{\tau_{2,y_{2}}^{-}%
<\tau_{2,b}^{+}}\right] \nonumber\\
&  =e^{\alpha b}\frac{W_{2}^{(q)}(x-y_{2})}{W_{2}^{(q)}(b-y_{2})}+e^{\alpha
y_{2}}\left(  Z_{2}^{(q)}(x-y_{2},\alpha)-\frac{W_{2}^{(q)}(x-y_{2})}%
{W_{2}^{(q)}(b-y_{2})}Z_{2}^{(q)}(b-y_{2},\alpha)\right)  .
\end{align}
Taking derivative of (\ref{eq.jointdef1}) with respect to $\alpha$ at
$\alpha=0,$ as in (53) of Avram et al. (2019),
\begin{align}
&  h_{2,2}(x,y_{2},b)=\frac{\partial}{\partial\alpha}\mathbb{E}_{x}%
^{2}[e^{\alpha X_{2,(\tau_{2,y_{2}}^{-}\wedge\tau_{2,b}^{+})}-q(\tau_{2,y_{2}%
}^{-}\wedge\tau_{2,b}^{+})}]|_{\alpha=0}=b\frac{W_{2}^{(q)}(x-y_{2})}%
{W_{2}^{(q)}(b-y_{2})}\nonumber\\
&  +y_{2}\left(  Z_{2}^{(q)}(x-y_{2})-\frac{W_{2}^{(q)}(x-y_{2})}{W_{2}%
^{(q)}(b-y_{2})}Z_{2}^{(q)}(b-y_{2})\right) \nonumber\\
&  +\overline{Z}_{2}^{(q)}(x-y_{2})-\varphi_{2}^{\prime}(0)\overline{W}%
_{2}^{(q)}(x-y_{2})\nonumber\\
&  -\frac{W_{2}^{(q)}(x-y_{2})}{W_{2}^{(q)}(b-y_{2})}\left(  \overline{Z}%
_{2}^{(q)}(b-y_{2})-\varphi_{2}^{\prime}(0)\overline{W}_{2}^{(q)}%
(b-y_{2})\right)  . \label{eq.h22}%
\end{align}
Combining (\ref{eq.H2}), (\ref{eq.h21}) and (\ref{eq.h22}), we have
\begin{equation}
H_{2}^{(q)}(x,y_{2},b)=\left(  a_{2}+\frac{c_{2}\varphi_{2}^{\prime}(0)}%
{q}\right)  h_{2,1}(x,y_{2},b)+\frac{c_{2}}{q}\left(  x-h_{2,2}(x,y_{2}%
,b)\right)  . \label{eq.H2final}%
\end{equation}

Next, we obtain $H_{1}^{(q)}(x,y_{1})$ for $0\leq x<y_{1}$ --the expected
discounted holding cost starting at inventory level $x$ at phase 1 until
reaching $y_{1}$:
\begin{equation}
H_{1}^{(q)}(x,y_{1})=a_{1}\mathbb{E}_{x}^{1}\left[  \int_{0}^{\kappa_{y_{1}%
}^{+}}e^{-qs}ds\right]  +c_{1}\mathbb{E}_{x}^{1}\left[  \int_{0}^{\kappa_{y_{1}}^{+}%
}e^{-qs}R_{s}ds\right]  . \label{eq.H1}%
\end{equation}
By (\ref{kypkappa}),
\begin{equation}
\mathbb{E}_{x}^{1}\left[  \int_{0}^{\kappa_{y_{1}}^{+}}e^{-qs}ds\right]  =\frac{1}%
{q}\left(  1-\frac{Z_{1}^{(q)}(x)}{Z_{1}^{(q)}(y_{1})}\right)  .
\label{eq.H1a1}%
\end{equation}
In order to obtain the second term on the right-hand side of (\ref{eq.H1}), we
apply the Kella-Whitt martingale (\ref{eq.kella-whitt}) for the process
$X_{1}$ with $\varphi_{1}(\alpha)=\log\mathbb{E}[e^{\alpha\left(
X_{1,1}-X_{1,0}\right)  }]$, $Y_{s}=L_{s}-qs/\alpha$ and $V_{s}=X_{1,s}%
+L_{s}-(q/\alpha)s=R_{s}-(q/\alpha)s$. Then
\begin{align}
&  \mathbb{E}_{x}^{1}\left[  (\varphi_{1}(\alpha)-q)\int_{0}^{\kappa_{y_{1}%
}^{+}}e^{\alpha R_{s}-qs}ds+e^{\alpha R_{0}}-e^{\alpha R_{\kappa_{y_{1}}^{+}%
}-q\kappa_{y_{1}}^{+}}\right. \nonumber\label{eq.kella1}\\
&  \left.  +\alpha\int_{0}^{\kappa_{y_{1}}^{+}}e^{\alpha R_{s}-qs}dL_{s}%
^{c}+\sum_{0\leq s\leq\kappa_{y_{1}}^{+}}e^{\alpha R_{s}-qs}(1-e^{-\alpha
\Delta L_{s}})\right]  =0.
\end{align}
Note that $R(\kappa_{y_{1}}^{+})=y_{1}$ and that $dL_{s}^{c}\neq0$ or $\Delta
L_{s}\neq0$ implies that $R_{s}=0$. Thus (\ref{eq.kella1}) reduces to%

\begin{align}
&  \mathbb{E}_{x}^{1}\left[  (\varphi_{1}(\alpha)-q)\int_{0}^{\kappa_{y_{1}%
}^{+}}e^{\alpha R_{s}-qs}ds+e^{\alpha x}-e^{\alpha y_{1}-q\kappa_{y_{1}}^{+}%
}\right. \nonumber\label{eq.kella11}\\
&  \left.  +\alpha\int_{0}^{\kappa_{y_{1}}^{+}}e^{-qs}dL_{s}^{c}+\sum_{0\leq
s\leq\kappa_{y_{1}}^{+}}e^{-qs}(1-e^{-\alpha\Delta L_{s}})\right]  =0.
\end{align}

By Eq. (80)-(81) in \cite{avram},
\begin{align}
&  \mathbb{E}_{x}^{1}\left[  \int_{0}^{\kappa_{y_{1}}^{+}}e^{-qs}dL_{s}\right]
\nonumber\\
&  =\frac{Z_{1}^{(q)}(x)}{Z_{1}^{(q)}(y_{1})}\left(  \overline{Z}_{1}%
^{(q)}(y_{1})+\varphi_{1}^{\prime}(0))/q\right)  -\left(  \overline{Z}%
_{1}^{(q)}(x)+\varphi_{1}^{\prime}(0))/q\right)  .
\end{align}

Taking derivative of (\ref{eq.kella11}) with respect to $\alpha$ at
$\alpha=0,$ and applying (\ref{kypkappa}) yields:%

\begin{align}
&  \varphi_{1}^{\prime}(0)\int_{0}^{\kappa_{y_{1}}^{+}}e^{-qs}ds+x-y_{1}%
\frac{Z_{1}^{(q)}(x)}{Z_{1}^{(q)}(y_{1})}\nonumber\label{eq.kella14}\\
&  \frac{Z_{1}^{(q)}(x)}{Z_{1}^{(q)}(y_{1})}\left(  \overline{Z}_{1}%
^{(q)}(y_{1})+\varphi_{1}^{\prime}(0)/q\right)  -\left(  \overline{Z}%
_{1}^{(q)}(x)+\varphi_{1}^{\prime}(0)/q\right) \nonumber\\
&  =q\mathbb{E}_{x}\left[  \int_{0}^{\kappa_{y_{1}}^{+}}R_{s}e^{-qs}ds\right]
.
\end{align}
Applying (\ref{kypkappa}) and after some algebra (\ref{eq.kella14}) yields:

%

\begin{equation}
\mathbb{E}_{x}\left[  \int_{0}^{\kappa_{y_{1}}^{+}}R_{s}e^{-qs}ds\right]
=\frac{Z_{1}^{(q)}(x)}{Z_{1}^{(q)}(y_{1})}\overline{\overline{W}}_{1}%
^{(q)}(y_{1})-\overline{\overline{W}}_{1}^{(q)}(x). \label{eq.H1c1}%
\end{equation}
By (\ref{eq.H1}), the expected discounted "fixed" part of the holding cost
until reaching $y_{1}$ is given by
\[
a_{1}\int_{0}^{\kappa_{y_{1}}^{+}}e^{-qt}dt=\frac{a_{1}}{q}\left(
1-\frac{Z_{1}^{(q)}(x)}{Z_{1}^{(q)}(y_{1})}\right)  .
\]
Applying (\ref{eq.H1}), (\ref{eq.H1a1}) and (\ref{eq.H1c1})
, we get
\begin{align}
&  H_{1}^{(q)}(x,y_{1})=\frac{a_{1}}{q}\left(  1-\frac{Z_{1}^{(q)}(x)}%
{Z_{1}^{(q)}(y_{1})}\right) \nonumber\label{eq.H1xy}\\
&  +c_{1}\left(  \frac{Z_{1}^{(q)}(x)}{Z_{1}^{(q)}(y_{1})}\overline
{\overline{W}}_{1}^{(q)}(y_{1})-\overline{\overline{W}}_{1}^{(q)}(x)\right)  .
\end{align}
In order to obtain $\mathcal{H}_{2}^{(q)}(x)$ --the expected discounted
holding cost starting at inventory level $x$ at phase $2$-, we first derive
$\mathbb{E}_{x}^{2}\left[  e^{-q\tau_{2,y_{2}}^{-}}1_{\tau_{2,y_{2}}^{-}%
<\tau_{b}^{+}}\mathcal{H}_{1}^{(q)}(X_{2,\tau_{2,y_{2}}^{-}})\right]  .$




For a function $g$ satisfying the conditions of Theorem 2 in Loeffen (2018),
let us define
\[
\Omega^{2}(g(x))=\mathbb{E}_{x}^{2}[e^{-q\tau_{2,y_{2}}^{-}}g(X_{2,\tau
_{2,y_{2}}^{-}})1_{\tau_{2,y_{2}}^{-}<\tau_{2,b}^{+}}].
\]
Then, by the aforementioned Theorem 2,
\begin{align}
&  \Omega^{2}(g(x))=g(x)-\frac{W_{2}^{(q)}(x-y_{2})}{W_{2}^{(q)}(b-y_{2}%
)}g(b)\nonumber\label{loeffen}\\
&  +\int_{y_{2}}^{b}(\mathcal{G}_{2}-q)g(z)\left(  \frac{W_{2}^{(q)}(x-y_{2}%
)}{W_{2}^{(q)}(b-y_{2})}W_{2}^{(q)}(b-z)-W_{2}^{(q)}(x-z)\right)  dz,
\end{align}
where $\mathcal{G}_{2}$ is the infinitesimal generator of $X_{2}$. Let us
first find $\Omega^{2}(Z_{1}^{q})(x)$: if $\mathcal{G}_{1}$ is the
infinitesimal generator of $X_{1}$, then
\begin{equation}
(\mathcal{G}_{2}-\mathcal{G}_{1})g(x)=(\sigma_{2}-\sigma_{1})g^{\prime}(x).
\label{eq.G2mG1}%
\end{equation}
It is well known that
\begin{equation}
(\mathcal{G}_{1}-q)Z_{1}^{(q)}(x)=0. \label{eq.GZ0}%
\end{equation}
Thus equations (\ref{loeffen}) and (\ref{eq.G2mG1}) yield:



%

\begin{align}
&  \Omega^{2}(Z_{1}^{(q)}(x))=Z_{1}^{(q)}(x)-\frac{W_{2}^{(q)}(x-y_{2})}%
{W_{2}^{(q)}(b-y_{2})}Z_{1}^{(q)}(b)\nonumber\\
&  +\int_{y_{2}}^{b}\left(  (\mathcal{G}_{1}-q)Z_{1}^{(q)}(z)+(\sigma
_{2}-\sigma_{1})qW_{1}^{(q)}(x)\right)  \left[  \frac{W_{2}^{(q)}(x-y_{2}%
)}{W_{2}^{(q)}(b-y_{2})}W_{2}^{(q)}(b-z)-W_{2}^{(q)}(x-z)\right]
dz\nonumber\\
&  =Z_{1}^{(q)}(x)-\frac{W_{2}^{(q)}(x-y_{2})}{W_{2}^{(q)}(b-y_{2})}%
Z_{1}^{(q)}(b)\nonumber\\
&  +\int_{y_{2}}^{b}(\sigma_{2}-\sigma_{1})qW_{1}^{(q)}(z)\left[  \frac
{W_{2}^{(q)}(x-y_{2})}{W_{2}^{(q)}(b-y_{2})}W_{2}^{(q)}(b-z)-W_{2}%
^{(q)}(x-z)\right]  dz. \label{eq.OmegaZ1}%
\end{align}

Similarly, let us consider
\begin{equation}
\Omega^{2}(\overline{\overline{W}}_{1}^{(q)}(x))=\mathbb{E}_{x}^{2}\left[
e^{-q\tau_{2,y_{2}}^{-}}1_{\tau_{2,y_{2}}^{-}<\tau_{2,b}^{+}}\overline
{\overline{W}}_{1}^{(q)}(X(\tau_{2,y_{2}}^{-}))\right]  ,
\end{equation}
by (\ref{eq.GZ0}) we have
\begin{equation}
(\mathcal{G}_{1}-q)\overline{\overline{W}}_{1}^{(q)}(x)=x.
\end{equation}
Thus,
\begin{align}
&  \Omega^{2}(\overline{\overline{W}}_{1}^{(q)}(x))=\overline{\overline{W}%
}_{1}^{(q)}(x)-\frac{W_{2}^{(q)}(x-y_{2})}{W_{2}^{(q)}(b-y_{2})}%
\overline{\overline{W}}_{1}^{(q)}(b)\nonumber\\
&  +\int_{y_{2}}^{b}(\mathcal{G}_{1}-q)\overline{\overline{W}}_{1}%
^{(q)}(z)+(\sigma_{2}-\sigma_{1})\overline{W}_{1}^{(q)}(z))\left[  \frac
{W_{2}^{(q)}(x-y_{2})}{W_{2}^{(q)}(b-y_{2})}W_{2}^{(q)}(b-z)-W_{2}%
^{(q)}(x-z)\right]  dz\nonumber\\
&  =\overline{\overline{W}}_{1}^{(q)}(x)-\frac{W_{2}^{(q)}(x-y_{2})}%
{W_{2}^{(q)}(b-y_{2})}\overline{\overline{W}}_{1}^{(q)}(b)\nonumber\\
&  +\int_{y_{2}}^{b}\left(  z+(\sigma_{2}-\sigma_{1})\overline{W}_{1}%
^{(q)}(z)\right)  \left[  \frac{W_{2}^{(q)}(x-y_{2})}{W_{2}^{(q)}(b-y_{2}%
)}W_{2}^{(q)}(b-z)-W_{2}^{(q)}(x-z)\right]  dz.\nonumber
\end{align}
Equations (\ref{eq.calH1}), (\ref{eq.H1xy}), (\ref{eq.calH2}) and
(\ref{eq.exitdown}) yield:%

\begin{align}
&  \mathcal{H}_{2}^{(q)}(x)=H_{2}(x,y_{2},b)\nonumber\label{eq.calH2main}\\
&  +\frac{a_{1}}{q}\left(  Z_{2}^{(q)}(x-y_{2})-\frac{W_{2}^{(q)}(x-y_{2}%
)}{W_{2}^{(q)}(b-y_{2})}Z_{2}^{(q)}(b-y_{2})-\frac{\Omega^{2}(Z_{1}^{(q)}%
(x))}{Z_{1}^{(q)}(y_{1})}\right) \nonumber\\
&  +c_{1}\left(  \frac{\Omega^{2}(Z_{1}^{(q)}(x))}{Z_{1}^{(q)}(y_{1}%
)}\overline{\overline{W}}_{1}^{(q)}(y_{1})-\Omega^{2}(\overline{\overline{W}%
}_{1}^{(q)}(x))\right) \nonumber\\
&  +\frac{\Omega^{2}(Z_{1}^{(q)}(x))}{Z_{1}^{(q)}(y_{1})}\mathcal{H}_{2}%
^{(q)}(y_{1})+\frac{W_{2}^{(q)}(x-y_{2})}{W_{2}^{(q)}(b-y_{2})}\mathcal{H}%
_{0}^{(q)}(b).
\end{align}
Let
\begin{align}
&  A(x):=H_{2}^{(q)}(x,y_{2},b)\nonumber\\
&  +\frac{a_{1}}{q}\left(  Z_{2}^{(q)}(x-y_{2})-\frac{W_{2}^{(q)}(x-y_{2}%
)}{W_{2}^{(q)}(b-y_{2})}Z_{2}^{(q)}(b-y_{2})-\frac{\Omega^{2}(Z_{1}^{(q)}%
(x))}{Z_{1}^{(q)}(y_{1})}\right) \nonumber\\
&  +c_{1}\left(  \frac{\Omega^{2}(Z_{1}^{(q)}(x))}{Z_{1}^{(q)}(y_{1}%
)}\overline{\overline{W}}_{1}^{(q)}(y_{1})-\Omega^{2}(\overline{\overline{W}%
}_{1}^{(q)}(x))\right)  ,
\end{align}
then
\begin{equation}
\mathcal{H}_{2}^{(q)}(x)=A(x)+\frac{\Omega^{2}(Z_{1}^{(q)}(x))}{Z_{1}%
^{(q)}(y_{1})}\mathcal{H}_{2}^{(q)}(y_{1})+\frac{W_{2}^{(q)}(x-y_{2})}%
{W_{2}^{(q)}(b-y_{2})}\mathcal{H}_{0}^{(q)}(b). \label{eq.calH2a}%
\end{equation}
Substituting $x=y_{1}$ in (\ref{eq.calH2a}) and solving for $\mathcal{H}%
_{2}^{(q)}(y_{1})$ yield:

\begin{equation}
\mathcal{H}_{2}^{(q)}(y_{1})=\frac{A(y_{1})+\frac{W_{2}^{(q)}(y_{1}-y_{2}%
)}{W_{2}^{(q)}(b-y_{2})}\mathcal{H}_{0}^{(q)}(b)}{1-\frac{\Omega^{2}%
(Z_{1}^{(q)}(y_{1}))}{Z_{1}^{(q)}(y_{1})}}. \label{eq.calHy1}%
\end{equation}
Equations (\ref{eq.calH2main})-(\ref{eq.calHy1}) yield that%

\begin{equation}
\mathcal{H}_{2}^{(q)}(x)=\alpha_{2}(x)\mathcal{H}_{0}^{(q)}(b)+\beta_{2}(x),
\label{eq.calH2linear}%
\end{equation}
where%

\begin{align}
&  \alpha_{2}(x)=\frac{W_{2}^{(q)}(x-y_{2})}{W_{2}^{(q)}(b-y_{2})}%
+\frac{\Omega^{2}(Z_{1}^{(q)}(x))}{Z_{1}^{(q)}(y_{1})-\Omega^{2}(Z_{1}%
^{(q)}(y_{1}))}\frac{W_{2}^{(q)}(y_{1}-y_{2})}{W_{2}^{(q)}(b-y_{2})}\\
& \nonumber\\
&  \beta_{2}(x)=A(x)+\frac{\Omega^{2}(Z_{1}^{(q)}(y_{1}))A(y_{1})}{Z_{1}%
^{(q)}(y_{1})-\Omega^{2}(Z_{1}^{(q)}(y_{1}))}.
\end{align}



Substituting (\ref{eq.calH2linear}) in (\ref{eq.calH1}), we get%

\begin{align}
&  \mathcal{H}_{1}^{(q)}(x)=H_{1}^{(q)}(x,y_{1})+\frac{Z_{1}^{(q)}(x)}%
{Z_{1}^{(q)}(y_{1})}\alpha_{2}(y_{1})\mathcal{H}_{0}^{(q)}(b)+\beta_{2}%
(y_{1})\frac{Z_{1}^{(q)}(x)}{Z_{1}^{(q)}(y_{1})}\nonumber\\
&  =\alpha_{1}(x)\mathcal{H}_{0}^{(q)}(b)+\beta_{1}(x),
\end{align}
where
\begin{align}
&  \alpha_{1}(x)=\frac{Z_{1}^{(q)}(x)}{Z_{1}^{(q)}(y_{1})}\alpha_{2}%
(y_{1})\nonumber\\
& \nonumber\\
&  \beta_{1}(x)=H_{1}^{(q)}(x,y_{1})+\frac{Z_{1}^{(q)}(x)}{Z_{1}^{(q)}(y_{1}%
)}\beta_{2}(y_{1}). \label{eq.calH1linear}%
\end{align}
In order to obtain $\mathcal{H}_{i}^{(q)}\ $for $i=1,2$, we substitute
(\ref{eq.calH1linear}) and (\ref{eq.calH2linear}) in (\ref{eq.calHb}) and get
the following linear equation for $\mathcal{H}_{0}^{(q)}(b)$.%

\begin{align}
&  \mathcal{H}_{0}^{(q)}(b)=\frac{h_{0}(b)}{q+\lambda}+\frac{\lambda}%
{\lambda+q}\left(  \int_{0}^{b-y_{3}}\mathcal{H}_{2}^{(q)}(b-z)dF(z)+\int
_{b-y_{3}}^{\infty}\mathcal{H}_{1}^{(q)}(b-z)dF(z)\right) \nonumber\\
&  =\frac{h_{0}(b)}{q+\lambda}+\frac{\lambda}{\lambda+q}\left(  \int
_{0}^{b-y_{3}}\beta_{2}(b-z)dF(z)+\int_{b-y_{3}}^{\infty}\beta_{1}%
(b-z)dF(z)\right) \nonumber\\
&  +\mathcal{H}_{0}^{(q)}(b)\left(  \int_{0}^{b-y_{3}}\alpha_{2}%
(b-z)dF(z)+\int_{b-y_{3}}^{\infty}\alpha_{1}(b-z)dF(z)\right)  .
\label{eq.Hbeq}%
\end{align}

We obtain $\mathcal{H}_{0}^{(q)}(b)$ solving the linear equation
(\ref{eq.Hbeq}); from this, we get $\mathcal{H}_{i}^{(q)}$ for $i=1,2$.

\subsubsection{Expected discounted shortage cost.}

Here, we derive formulas for the expected discounted shortage cost
$\mathcal{S}_{i}^{(q)}(x)$ starting at inventory level $x$ at phase $i$ for
$i=1,2$ together with the expected discounted shortage cost $\mathcal{S}%
_{0}^{(q)}(b)$ starting at inventory level $b$.

Let us define $S_{1}^{(q)}(x,y_{1})$ as the expected discounted shortage cost
starting at inventory level $x\in\lbrack0,y_{1})$ until the inventory level
reaches $y_{1}$.
By (\ref{kypkappa}), we can write%

\begin{equation}
\mathcal{S}_{1}^{(q)}(x)=S_{1}^{(q)}(x,y_{1})+\mathbb{E}_{x}^{1}%
[e^{-\kappa_{y_{1}}^{+}}]\mathcal{S}_{2}^{(q)}(y_{1})=S_{1}^{(q)}%
(x,y_{1})+\frac{Z_{1}^{(q)}(x)}{Z_{1}^{(q)}(y_{1})}\mathcal{S}_{2}^{(q)}%
(y_{1}). \label{eq.S1cal}%
\end{equation}
Equations (\ref{eq.u}) and (\ref{eq.exitup}) yield
\begin{align}
&  \mathcal{S}_{2}^{(q)}(x)=\int_{y_{2}}^{b}u_{2}^{(q)}(y_{2},b,x,z)\lambda
\left(  \int_{v=z}^{\infty}\left(  p(v-z)+\mathcal{S}_{1}^{(q)}(0)\right)
dF(v)\right)  dz\nonumber\label{eq.S2cal}\\
&  +\int_{y_{2}}^{b}u_{2}^{(q)}(y_{2},b,x,z)\lambda\left(  \int_{v=z-y_{2}%
}^{z}\mathcal{S}_{1}^{(q)}(z-v)dF(v)\right)  dz\nonumber\\
&  +\frac{W_{2}^{(q)}(x-y_{2})}{W_{2}^{(q)}(b-y_{2})}\mathcal{S}_{0}^{(q)}(b).
\end{align}
And also,
\begin{align}
&  \mathcal{S}_{0}^{(q)}(b)=\frac{\lambda}{\lambda+q}\left(  \int_{b}^{\infty
}(p(z-b)+\mathcal{S}_{1}^{(q)}(0))dF(z)\right. \nonumber\\
&  +\left.  \int_{0}^{b-y_{3}}\mathcal{S}_{2}^{(q)}(b-z)dF(z)+\int_{b-y_{3}%
}^{b}\mathcal{S}_{1}^{(q)}(b-z)dF(z)\right)  . \label{eq.calSb}%
\end{align}

First we consider $S_{1}^{(q)}(x,y_{1})$ which corresponds to the expected
discounted shortage cost starting at $x$ at phase 1, $0<x<y_{1}$ until the
process reaches $y_{1}$. The definition of $u_{1}^{(q)}$ given in (\ref{eq.u})
yields
%

\begin{align}
&  S_{1}^{(q)}(x,y_{1})=\int_{0}^{y_{1}}u_{1}^{(q)}(0,y_{1},x,z)\lambda\left(
\int_{v=z}^{\infty}p(v-z)dF(v)\right)  dz\label{eq.S1xy1}\\
&  +\mathbb{E}_{x}^{1}\left(  e^{-q\tau_{1,0}^{-}}1_{\tau_{1,0}^{-}%
<\tau_{1,y_{1}}^{+}}\right)  \cdot\int_{0}^{y_{1}}u_{1}^{(q)}(0,y_{1}%
,0,z)\lambda\left(  \int_{v=z}^{\infty}p(v-z)dF(v)\right)  dz \label{eq.S1xy2}%
\\
&  \cdot\sum_{j=0}^{\infty}\left(  \mathbb{E}_{0}^{1}[e^{-q\tau_{1,0}^{-}%
}1_{\tau_{1,0}^{-}<\tau_{1,y_{1}}^{+}}]\right)  ^{j},\nonumber
\end{align}
where (\ref{eq.S1xy1}) describes the expected discounted shortage cost
occurring before the inventory level reaches $y_{1}$, and (\ref{eq.S1xy2})
describes the expected discounted shortage costs occurring after the first
downcrossing level $0$. Applying equation (\ref{eq.exitdown}) yields:
\begin{align}
&  S_{1}^{(q)}(x,y_{1})=\int_{0}^{y_{1}}u_{1}^{(q)}(0,y_{1},x,z)\lambda\left(
\int_{v=z}^{\infty}p(v-z)dF(v)\right)  dz\nonumber\\
&  +\frac{Z^{(q)}(x)-\frac{W_{1}^{(q)}(x)}{W_{1}^{(q)}(y_{1})}Z_{1}%
^{(q)}(y_{1})}{\frac{W_{1}^{(q)}(0)}{W_{1}^{(q)}(y_{1})}Z_{1}^{(q)}(y_{1}%
)}\int_{0}^{y_{1}}u_{1}^{(q)}(0,y_{1},0,z)\lambda\left(  \int_{v=z}^{\infty
}p(v-z)dF(v)\right)  dz.
\end{align}




Next, we obtain linear equations to obtain $\mathcal{S}_{i}^{(q)}(x)$ and
$\mathcal{S}_{0}^{(q)}(b)$. Let us define
\begin{align}
&  \mu(x):=\int_{y_{2}}^{b}u_{2}^{(q)}(y_{2},b,x,z)\lambda\left(  \int
_{v=z}^{\infty}p(v-z)dF(z)\right)  dz\nonumber\\
&  +\lambda S_{1}^{(q)}(0,y_{1})\int_{y_{2}}^{b}u_{2}^{(q)}(y_{2}%
,b,x,z)\overline{F}(z)dz\nonumber\\
&  +\int_{y_{2}}^{b}u_{2}^{(q)}(y_{2},b,x,z)\lambda\left(  \int_{v=z-y_{2}%
}^{z}S_{1}^{(q)}(z-v,y_{1})dF(v)\right)  dz,
\end{align}
where $\bar{F}(z)=1-F(z)$. Let us also define
\begin{align}
&  \gamma(x):=\lambda\frac{1}{Z_{1}^{(q)}(y_{1})}\int_{y_{2}}^{b}u_{2}%
^{(q)}(y_{2},b,x,z)\overline{F}(z)dz\nonumber\label{eq.gamma}\\
&  +\lambda\int_{y_{2}}^{b}u_{2}^{(q)}(y_{2},b,x,z)\left(  \int_{v=z-y_{2}%
}^{z}\frac{Z_{1}^{(q)}(z-v)}{Z_{1}^{(q)}(y_{1})}dF(z)\right)  dz.
\end{align}
Substituting (\ref{eq.S1cal}) in (\ref{eq.S2cal}), we have that $\mathcal{S}%
_{2}^{\left(  q\right)  }(x)$ can be written as follows:%

\begin{equation}
\mathcal{S}_{2}^{\left(  q\right)  }(x)=\mu(x)+\gamma(x)\mathcal{S}%
_{2}^{\left(  q\right)  }(y_{1})+\frac{W_{2}^{(q)}(x-y_{2})}{W_{2}%
^{(q)}(b-y_{2})}\mathcal{S}_{0}^{\left(  q\right)  }(b). \label{eq.calS2lin}%
\end{equation}
Solving (\ref{eq.calS2lin}) for $x=y_{1}$ yield
\begin{equation}
\mathcal{S}_{2}^{\left(  q\right)  }(y_{1})=\frac{\mu(y_{1})+\frac{W_{2}%
^{(q)}(y_{1}-y_{2})}{W_{2}^{(q)}(b-y_{2})}\mathcal{S}_{0}^{\left(  q\right)
}(b)}{1-\gamma(y_{1})}. \label{eq.calS2y1}%
\end{equation}
Let us define
\[
\mu_{2}(x):=\mu(x)+\gamma(x)\frac{\mu(y_{1})}{1-\gamma(y_{1})},
\]

\[
\gamma_{2}(x):=\frac{W_{2}^{(q)}(x-y_{2})}{W_{2}^{(q)}(b-y_{2})}+\frac
{W_{2}^{(q)}(y_{1}-y_{2})}{W_{2}^{(q)}(b-y_{2})}\frac{\gamma(x)}%
{1-\gamma(y_{1})},
\]
and
\begin{align}
&  \mu_{1}(x):=S_{1}^{\left(  q\right)  }(x,y_{1})+\frac{Z_{1}^{(q)}(x)}%
{Z_{1}^{(q)}(y_{1})}\mu_{2}(y_{1}),\nonumber\\
&  \gamma_{1}(x):=\frac{Z_{1}^{(q)}(x)}{Z_{1}^{(q)}(y_{1})}\gamma_{2}(y_{1}).
\end{align}
Thus, equations (\ref{eq.S2cal}), (\ref{eq.S1cal}), (\ref{eq.calS2lin}) and
(\ref{eq.calS2y1}) yield:
\begin{equation}
\mathcal{S}_{i}^{\left(  q\right)  }(x)=\mu_{i}(x)+\mathcal{S}_{0}^{\left(
q\right)  }(b)\gamma_{i}(x)\ \text{for}\,\,\,i=1,2. \label{eq.calSi}%
\end{equation}
Substituting (\ref{eq.calSi}) in (\ref{eq.calSb}), we get the following linear
equation for $\mathcal{S}_{0}(b)$,
\begin{align}
&  \mathcal{S}_{0}^{(q)}(b)=\frac{\lambda}{\lambda+q}\left(  \int_{b}^{\infty
}(p(z-b)+\mu_{1}(0)+\gamma_{1}(0)\mathcal{S}_{0}^{(q)}(b))dF(z)\right.
\nonumber\\
&  +\left.  \int_{0}^{b-y_{3}}(\mu_{2}(b-z)+\mathcal{S}_{0}^{(q)}(b)\gamma
_{2}(b-z))dF(z)+\int_{b-y_{3}}^{b}(\mu_{1}(b-z)+\gamma_{1}(b-z)\mathcal{S}%
_{0}^{(q)}(b))dF(z)\right)  ,\nonumber
\end{align}
and so
\begin{equation}
\mathcal{S}_{0}^{(q)}(b)=\frac{\frac{\lambda}{\lambda+q}\left(  \int
_{b}^{\infty}(p(z-b)+\mu_{1}(0))dF(z)+\int_{0}^{b-y_{3}}\mu_{2}(b-z)dF(z)+\int
_{b-y_{3}}^{b}\mu_{1}(b-z)dF(z)\right)  }{1-\gamma_{1}(0)\overline{F}%
(b)+\int_{0}^{b-y_{3}}\gamma_{2}(b-z)dF(z)+\int_{b-y_{3}}^{b}\gamma
_{1}(b-z)dF(z)}.
\end{equation}

Finally, from (\ref{eq.calSi}), we get the formulas for $\mathcal{S}%
_{i}^{\left(  q\right)  }(x)\ $for$\,\,\,i=1,2$.

\subsubsection{Expected discounted switching cost.}

Here, we compute the formulas for the expected discounted switching cost. Let
$\mathcal{K}_{i}^{(q)}(x)$, $i=1,2$ be the expected discounted switching cost
starting at inventory level $x$ and phase $i$ for $i=1,2$ and let
$\mathcal{K}_{0}^{(q)}(b)$ be the expected discounted switching cost starting
at $b$.
Assume that initially the inventory level $x$ is at phase $1$, then the first
switching from phase $1$ to phase $2$ occurs at $\kappa_{y_{1}}^{+}$. By
(\ref{kypkappa}),
\begin{equation}
\mathcal{K}_{1}^{(q)}(x)=\mathbb{E}_{x}^{1}[e^{-a\kappa_{y_{1}}^{+}}]\left(
K_{12}+\mathcal{K}_{2}^{(q)}(y_{1})\right)  =\frac{Z_{1}^{(q)}(x)}{Z_{1}%
^{(q)}(y_{1})}\left(  K_{12}+\mathcal{K}_{2}^{(q)}(y_{1})\right)  .
\label{eq.calK1}%
\end{equation}

If initially the inventory level $x$ is at phase $2$, then the first switching
from phase $2$ to phase $1$ occurs when the inventory level downcrosses
$y_{2}$ before reaching $b$. If the inventory reaches $b$ before downcrossing
$y_{2}$, there is a switching from phase $2$ to phase $0$. By (\ref{eq.exitup}%
),$\ $
\begin{equation}
\mathcal{K}_{2}^{(q)}(x)=\frac{W_{2}^{(q)}(x-y_{2})}{W_{2}^{(q)}(b-y_{2}%
)}\left(  K_{20}+\mathcal{K}_{0}^{(q)}(b)\right)  +\mathbb{E}_{x}^{2}\left[
e^{-q\tau_{2,y_{2}}^{-}}1_{\tau_{2,y_{2}}^{-}<\tau_{2,b}^{+}}\,(K_{21}%
+\mathcal{K}_{1}^{(q)}(X_{2,\tau_{2,y_{2}}^{-}}))\right]  .
\end{equation}
Due to equations (\ref{eq.exitdown}), (\ref{eq.calK1}) and (\ref{eq.OmegaZ1}%
),
\begin{align}
&  \mathcal{K}_{2}^{(q)}(x)=\frac{W_{2}^{(q)}(x-y_{2})}{W_{2}^{(q)}(b-y_{2}%
)}\left(  K_{20}+\mathcal{K}_{0}^{(q)}(b)\right)  +\frac{\Omega^{2}%
(Z_{1}^{(q)}(x))}{Z_{1}^{(q)}(y_{1})}\left(  K_{12}+\mathcal{K}_{2}%
^{(q)}(y_{1})\right) \nonumber\label{eq.calK22}\\
&  +K_{21}\left(  Z_{2}^{(q)}(x-y_{2})-\frac{W_{2}^{(q)}(x-y_{2})}{W_{2}%
^{(q)}(b-y_{2})}Z_{2}^{(q)}(b-y_{2})\right)  .
\end{align}
Substituting $x\ $by $y_{1}$ in (\ref{eq.calK22}) and solving for
$\mathcal{K}_{2}^{(q)}(y_{1})$ yields:%

\begin{align}
&  \mathcal{K}_{2}^{(q)}(y_{1})=\frac{1}{1-\frac{\Omega^{2}(Z_{1}^{(q)}%
(y_{1}))}{Z_{1}^{(q)}(y_{1})}}\left(  K_{20}\frac{W_{2}^{(q)}(y_{1}-y_{2}%
)}{W_{2}^{(q)}(b-y_{2})}+K_{12}\frac{\Omega^{2}(Z_{1}^{(q)}(y_{1}))}%
{Z_{1}^{(q)}(y_{1})}\right. \nonumber\\
&  +K_{21}(Z_{2}^{(q)}(y_{1}-y_{2})-\frac{W_{2}^{(q)}(y_{1}-y_{2})}%
{W_{2}^{(q)}(b-y_{2})}Z_{2}^{(q)}(b-y_{2}))\nonumber\\
&  +\left.  \frac{W_{2}^{(q)}(y_{1}-y_{2})}{W_{2}^{(q)}(b-y_{2})}%
\mathcal{K}_{0}^{(q)}(b)\right)  .
\end{align}
Thus
\begin{equation}
\mathcal{K}_{2}^{(q)}(x)=\omega_{2}(x)+\delta_{2}(x)\mathcal{K}_{0}^{(q)}(b),
\end{equation}
where,%

\begin{align}
&  \omega_{2}(x):=K_{20}\frac{W_{2}^{(q)}(x-y_{2})}{W_{2}^{(q)}(b-y_{2}%
)}+K_{21}\left(  Z_{2}^{(q)}(x-y_{2})-\frac{W^{(q)}(x-y_{2})}{W^{(q)}%
(b-y_{2})}Z_{2}^{(q)})(b-y_{2})\right) \nonumber\\
&  +\frac{\Omega^{2}(Z_{1}^{(q)}(x))}{Z_{1}^{(q)}(y_{1})}\left( K_{12}\right.
\nonumber\\
&  +\left.    \frac{K_{20}\frac{W_{2}^{(q)}(y_{1}-y_{2})}{W_{2}%
^{(q)}(b-y_{2})}+K_{12}\frac{\Omega^{2}(Z_{1}^{(q)}(y_{1}))}{Z_{1}^{(q)}%
(y_{1})}+K_{21}\left(Z_{2}^{(q)}(y_{1}-y_{2})-\frac{W_{2}^{(q)}(y_{1}-y_{2})}%
{W_{2}^{(q)}(b-y_{2})}Z_{2}^{(q)}(b-y_{2})\right)}{1-\frac{\Omega^{2}(Z_{1}%
^{(q)}(y_{1}))}{Z_{1}^{(q)}(y_{1})}}\right)  
\end{align}
and
\begin{equation}
\delta_{2}(x):=\frac{W_{2}^{(q)}(x-y_{2})}{W_{2}^{(q)}(b-y_{2})}+\frac
{\frac{\Omega^{2}(Z_{1}^{(q)}(x))}{Z_{1}^{(q)}(y_{1})}\frac{W_{2}^{(q)}%
(y_{1}-y_{2})}{W_{2}^{(q)}(b-y_{2})}}{1-\frac{\Omega^{2}(Z_{1}^{(q)}(y_{1}%
))}{Z_{1}^{(q)}(y_{1})}}.
\end{equation}
By (\ref{eq.calK1}),
\begin{equation}
\mathcal{K}_{1}^{(q)}(x)=\omega_{1}(x)+\delta_{1}(x)\mathcal{K}_{0}^{(q)}(b),
\end{equation}
where
\[
\omega_{1}(x):=\frac{Z_{1}^{(q)}(x)}{Z_{1}^{(q)}(y_{1})}\left(  K_{12}%
+\omega_{2}(y_{1})\right)
\]
and
\begin{equation}
\delta_{1}(x):=\frac{Z_{1}^{(q)}(x)}{Z_{1}^{(q)}(y_{1})}\delta_{2}(y_{1}).
\end{equation}

Moreover, $\mathcal{K}_{0}^{(q)}(b)$ satisfies the following linear equation:%

\begin{align}
&  \mathcal{K}_{0}^{(q)}(b)=\frac{\lambda}{q+\lambda}\left(  \int_{0}%
^{b-y_{3}}(K_{02}+\mathcal{K}_{2}^{(q)}(b-z))dF(z)\right. \nonumber\\
&  \left.  +\int_{b-y_{3}}^{b}(K_{01}+\mathcal{K}_{1}^{(q)}(b-z))dF(z)+\int
_{b}^{\infty}(K_{01}+\mathcal{K}_{1}^{(q)}(0))dF(z)\right) \nonumber\\
&  =\frac{\lambda}{q+\lambda}\left(  \int_{0}^{b-y_{3}}(K_{02}+\omega
_{2}(b-z)+\delta_{2}(b-z)\mathcal{K}_{0}^{(q)}(b))dF(z)\right. \nonumber\\
&  +\int_{b-y_{3}}^{b}(K_{01}+\omega_{1}(b-z)+\delta_{1}(b-z)\mathcal{K}%
_{0}^{(q)}(b))dF(z)\nonumber\\
&  +\left.  \int_{b}^{\infty}(K_{01}+\omega_{1}(0)+\delta_{1}(0)\mathcal{K}%
_{0}^{(q)}(b))dF(z)\right)  ,
\end{align}
thus
\begin{equation}
\mathcal{K}_{0}^{(q)}(b)=\frac{\frac{\lambda}{\lambda+q}\left(  \int
_{0}^{b-y_{3}}(K_{02}+\omega_{2}(b-z))dF(z)+\int_{b-y_{3}}^{b}(K_{01}%
+\omega_{1}(b-z))dF(z)+(K_{01}+\omega_{1}(0))\overline{F}(b)\right)  }%
{1-\frac{\lambda}{\lambda+q}\left(  \int_{0}^{b-y_{3}}\delta_{2}%
(b-z)dF(z)+\int_{b-y_{3}}^{b}\delta_{1}(b-z)dF(z)+\delta_{1}(0)\overline
{F}(b)\right)  .}.
\end{equation}

\subsubsection{Total discounted cost.}

As a summary, we have that the total discounted cost starting at inventory
level $x$ and phase $i$ is: \begin{table}[th]
\centering	
\begin{tabular}
[c]{ccc}%
Inventory level & Phase & Expected discounted cost\\\hline
$0\leq x<y_{1}$ & 1 & $\mathcal{H}_{1}^{(q)}(x)+\mathcal{S}_{1}^{(q)}%
(x)+\mathcal{K}_{1}^{(q)}(x)$\\
$0\leq x\leq y_{2}$ & 2 & $\mathcal{H}_{1}^{(q)}(x)+\mathcal{S}_{1}%
^{(q)}(x)+\mathcal{K}_{1}^{(q)}(x)+K_{21}$\\
$y_{2}<x<b$ & 2 & $\mathcal{H}_{2}^{(q)}(x)+\mathcal{S}_{2}^{(q)}%
(x)+\mathcal{K}_{2}^{(q)}(x)$\\
$y_{1}\leq x<b$ & 1 & $\mathcal{H}_{2}^{(q)}(x)+\mathcal{S}_{2}^{(q)}%
(x)+\mathcal{K}_{2}^{(q)}(x)+K_{12}$\\
$x=b$ & 0 & $\mathcal{H}_{0}^{(q)}(b)+\mathcal{S}_{0}^{(q)}(b)+\mathcal{K}%
_{0}^{(q)}(b)$\\
&  &
\end{tabular}
\end{table}

\subsection{Cost functions for strategies of type two.}

Here the switching zone from $1$ to $2$ is $A_{12}=[y_{1},y_{4}]$, the
switching zone from $2$ to $1$ is $A_{21}=[0,y_{2}]$, the selection zones are
$C_{1}=[0,y_{3}]$ and $C_{2}=(y_{3},b)$ and the non-action zone $(y_{2}%
,y_{1})\cup(y_{4},b)$ for $0\leq y_{2}\leq y_{3}<y_{1}<y_{4}<b$.

The analysis of the value function in this case is very similar to the
analysis of the strategy of type one. The only difference is in the case when
initially the inventory level is $x\in(y_{4},b)$ at phase $1$. Thus we
consider only this case.

Let us start with $\bar{\mathcal{H}}_{1}^{(q)}(x)$ ---the expected discounted
holding cost in the case $y_{4}<x<b$. Consider the expected discounted holding
until reaching $b$ or down-crossing $y_{4},$
\begin{equation}
H_{1}^{(q)}(x,y_{4},b)=\mathbb{E}_{x}^{1}\left[  \int_{0}^{\tau_{1,y_{4}}%
^{-}\wedge\tau_{1,b}^{+}}e^{-qs}(a_{1}+c_{1}X_{1,s})ds\right]  .
\end{equation}
Similarly to equations (\ref{eq.H2})-(\ref{eq.H2final}), we have
\begin{equation}
H_{1}^{(q)}(x,y_{4},b)=\left(  a_{1}+\frac{c_{1}\varphi_{1}^{\prime}(0)}%
{q}\right)  h_{1,1}(x,y_{4},b)+\frac{c_{1}}{q}\left(  x-h_{1,2}(x,y_{4}%
,b)\right)  ,
\end{equation}
where%

\begin{align}
&  h_{1,1}(x,y_{4},b)=\frac{1}{q}\left(  1-\mathbb{E}_{x}^{1}\left[
e^{-q\tau_{1,y_{4}}^{-}}\,1_{\tau_{1,y_{4}^{-}}<\tau_{1.b}^{+}}\right]
-\mathbb{E}_{x}^{1}\left[  e^{-q\tau_{1,b}^{+}}\,1_{\tau_{1.b}^{+}%
<\tau_{1,y_{4}^{-}}}\right]  \right) \nonumber\\
&  =\frac{1}{q}\left(  1-Z_{1}^{(q)}(x-y_{4})+\frac{W_{1}^{(q)}(x-y_{4}%
)}{W_{1}^{(q)}(b-y_{4})}Z_{1}^{(q)}(b-y_{4})-\frac{W_{1}^{(q)}(x-y_{4})}%
{W_{1}^{(q)}(b-y_{4})}\right)  \label{eq.h11}%
\end{align}
and
\begin{align}
&  h_{1,2}(x,y_{4},b)=\frac{\partial}{\partial\alpha}\mathbb{E}_{x}\left[
e^{\alpha(X_{1,\tau_{1,y_{4}}^{-}\wedge\tau_{1,b}^{+}})-q(\tau_{1,y_{4}}%
^{-}\wedge\tau_{1,b}^{+})}\right]  |_{\alpha=0}=b\frac{W_{1}^{(q)}(x-y_{4}%
)}{W_{1}^{(q)}(b-y_{4})}\nonumber\\
&  +y_{4}\left(  Z_{1}^{(q)}(x-y_{4})-\frac{W_{1}^{(q)}(x-y_{4})}{W_{1}%
^{(q)}(b-y_{4})}Z_{1}^{(q)}(b-y_{4})\right) \nonumber\\
&  +\overline{Z}_{1}^{(q)}(x-y_{4})-\varphi_{1}^{\prime}(0)\overline{W}%
_{1}^{(q)}(x-y_{4})\nonumber\\
&  -\frac{W_{1}^{(q)}(x-y_{4})}{W_{1}^{(q)}(b-y_{4})}\left(  \overline{Z}%
_{1}^{(q)}(b-y_{4})-\varphi_{1}^{\prime}(0)\overline{W}_{1}^{(q)}%
(b-y_{4})\right)  . \label{eq.h12}%
\end{align}
Once the inventory level reaches $b$, the expected discounted holding cost is
$\mathcal{H}_{0}^{(q)}(b)$. In the case that the inventory level down-crosses
$y_{4}$ before reaching $b$ there are two scenarios: \textbf{1.} If
$X_{1,\tau_{1,y_{4}}^{-}}$ lies in $[y_{1},y_{4}]$, then the expected
discounted holding cost is $\mathcal{H}_{2}^{(q)}\left(  X_{1,\tau_{1,y_{4}%
}^{-}}\right)  $. \textbf{2.} If the inventory level immediately after the
jump is $X_{1,\tau_{1,y_{4}}^{-}}$ lies in $(-\infty,y_{1})$, then the
expected discounted holding cost is $\mathcal{H}_{1}^{(q)}(X_{1,\tau_{1,y_{4}%
}^{-}})$. Thus,
\begin{align}
&  \bar{\mathcal{H}}_{1}^{(q)}(x)=H_{1}^{(q)}(x,y_{4},b)+\frac{W_{1}%
^{(q)}(x-y_{4})}{W_{1}^{(q)}(b-y_{4})}\cdot\mathcal{H}_{0}^{(q)}(b)\nonumber\\
&  +\lambda\int_{y_{4}}^{b}u_{1}^{(q)}(y_{4},b,x,y)\left(  \int_{z=y-y_{4}%
}^{y-y_{1}}\mathcal{H}_{2}^{(q)}(y-z)dF(z)\right)  dy\nonumber\\
&  +\lambda\int_{y_{4}}^{b}u_{1}^{(q)}(y_{4},b,x,y)\left(  \int_{z=y-y_{1}%
}^{y}\mathcal{H}_{1}^{(q)}(y-z)dF(z)\right)  dy\nonumber\\
&  +\mathcal{H}_{1}^{(q)}(0)\lambda\int_{y_{4}}^{b}u_{1}^{(q)}(y_{4}%
,b,x,y)\bar{F}(y)dy.
\end{align}

Consider now $\bar{\mathcal{S}}_{1}^{(q)}(x)$ ---the expected discounted
shortage cost starting at $x\in(y_{4},b)$ at phase $1$. If the process reaches
$b$ before down-crossing $y_{4}$, then the expected discounted shortage cost
is $\mathcal{S}_{0}^{(q)}(b)$. If the inventory level down-crosses $y_{4}$
before reaching $b$, then the shortage cost is $\mathcal{S}_{2}^{(q)}%
(X_{1,\tau_{1,y_{4}}^{-}})$ in the case that $y_{1}\leq X_{1,\tau_{1,y_{4}%
}^{-}}\leq y_{4}$, is $\mathcal{S}_{1}^{(q)}(X_{1,\tau_{1,y_{4}}^{-}})$ in the
case that $0\leq X_{1,\tau_{1,y_{4}}^{-}}<y_{1},$ and is $p(-X_{1,\tau
_{1,y_{4}}^{-}})+\mathcal{S}_{1}^{(q)}(0)$ in the case that $X_{1,\tau
_{1,y_{4}}^{-}}<0$. Applying (\ref{eq.exitup}) and (\ref{eq.u}) yields:
\begin{align}
&  \bar{\mathcal{S}}_{1}^{(q)}(x)=\frac{W_{1}^{(q)}(x-y_{4})}{W_{1}%
^{(q)}(b-y_{4})}\mathcal{S}_{0}^{(q)}(b)\nonumber\\
&  +\lambda\int_{y_{4}}^{b}u_{1}^{(q)}(y_{4},b,x,y)\int_{z=y-y_{4}}^{y-y_{1}%
}\mathcal{S}_{2}^{(q)}(y-z)dF(z)dy\nonumber\\
&  +\lambda\int_{y_{4}}^{b}u_{1}^{(q)}(y_{4},b,x,y)\int_{z=y-y_{1}}%
^{y}\mathcal{S}_{1}^{(q)}(y-z)dF(z)dy\nonumber\\
&  +\lambda\int_{y_{4}}^{b}u_{1}^{(q)}(y_{4},b,x,y)\int_{z=y}^{\infty
}p(z-y)dF(z)dy\nonumber\\
&  +\mathcal{S}_{1}^{(q)}(0)\lambda\int_{y_{4}}^{b}u_{1}^{(q)}(y_{4}%
,b,x,y)\bar{F}(y)dy.
\end{align}

Similarly we obtain the expected discounted switching cost $\bar{\mathcal{K}%
}_{1}^{(q)}(x)$:%

\begin{align}
&  \bar{\mathcal{K}}_{1}^{(q)}(x)=\frac{W_{1}^{(q)}(x-y_{4})}{W_{1}%
^{(q)}(b-y_{4})}\mathcal{K}_{0}^{(q)}(b)\nonumber\\
&  +\lambda\int_{y_{4}}^{b}u_{1}^{(q)}(x-y_{4},y-y_{4})\left(  \int
_{z=y-y_{4}}^{y-y_{1}}\mathcal{K}_{2}^{(q)}(y-z)dF(z)\right)  dy\nonumber\\
&  +K_{12}\lambda\int_{y_{4}}^{b}u_{1}^{(q)}(y_{4},b,x,y)(F(y-y_{4}%
)-F(y-y_{1}))dy\nonumber\\
&  +\lambda\int_{y_{4}}^{b}u_{1}^{(q)}(y_{4},b,x,y)\left(  \int_{z=y-y_{1}%
}^{y}\mathcal{K}_{1}^{(q)}(y-z)dF(z)\right)  dy\nonumber\\
&  +\mathcal{K}_{1}^{(q)}(0)\lambda\int_{y_{4}}^{b}u_{1}^{(q)}(y_{4}%
,b,x,y)\bar{F}(y)dy.
\end{align}

\section{Examples\label{Section Examples}}

In this section, we find the optimal strategies for three different
situations. In the first one, the optimal strategy is of Doshi type, in the
second is of type one and in the third is of type two.

\subsection{First Example: Doshi strategy is optimal}

In this example we consider two equal manufacturing units, in phase $1$ both
units are producing together, in phase $2$ only one manufacturing unit is
working and in phase $0$ none of the units are working. We assume that the
cost of shutdown each unit is equal to $1$, the cost of restarting each unit
is equal to $2$ and the rate of production of each unit is equal to $3/2$. The
cost rate of production of each unit is $1/100,$ and there is also a fix cost
rate (independent of the production) equal to $1/1000$. Moreover, the holding
cost rate is $1/1000$ and the storage capacity is $b=10$. So we have the
following parameters $\sigma_{1}=3$, $\sigma_{2}=3/2$, $K_{12}=1$,~$K_{21}=2$,
$K_{20}=2$, $K_{10}=4$,$~K_{02}=2,$ $K_{01}=4$, $h_{2}(x)=(21+x)/1000$,
$h_{1}(x)=(41+x)/1000$, $h_{0}(b)=(1+b)/1000$. We assume that the rate of
arrival of the customer demands is $\lambda=2$, the demands are distributed as
$Exp(1.5)$, the discount rate is$\ q=0.1$ and the penalty cost when an amount
$y$ of a costumer is lost is given by $p(y)=(80+40y)/100$ (here we are taking
$l=0)$.

We find that the best Doshi strategy is given by the sets $A_{{\small 12}%
}=[y_{{\small 1}},10)$, $A_{{\small 21}}=[0,y_{{\small 2}}]$, $C_{{\small 1}%
}=[0,y_{{\small 2}}]$ and $C_{{\small 2}}=(y_{{\small 2}},10)$ with
$y_{{\small 2}}=1.526\ $and $y_{{\small 1}}=5.077$. We check that the value
functions of this strategy are viscosity solutions of the equations
(\ref{HJBi}) and satisfy the conditions of Theorem \ref{Teorema Verificacion};
so the best Doshi strategy is the optimal one. We show this optimal strategy
en Figure 1.%

\[%
{\parbox[b]{4.8512in}{\begin{center}
\includegraphics[
height=2.419in,
width=4.8512in
]%
{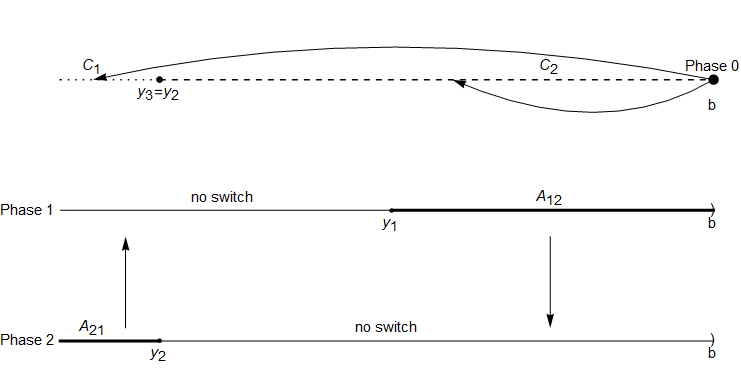}%
\\
Figure 1: Optimal strategy in first example.
\end{center}}}%
\]

In Figure 2, we show the discounted total cost $V_{1}(x)$ (dotted), $V_{2}(x)$
(dashed) and $V_{0}(b)$ (solid point) of the optimal strategy; in Figure 3,
the discounted holding cost $\mathcal{H}_{1}^{(q)}(x)$ (dotted),
$\mathcal{H}_{2}^{(q)}(x)$ (dashed) and $\mathcal{H}_{0}^{(q)}(b)$ (solid
point) of the optimal strategy; in Figure 4, the discounted penalty cost
$\mathcal{S}_{1}^{(q)}(x)$ (dotted), $\mathcal{S}_{2}^{(q)}(x)$ (dashed) and
$\mathcal{S}_{0}^{(q)}(b)$ (solid point) of the optimal strategy; and finally
in Figure 5, the discounted penalty cost $\mathcal{K}_{1}^{(q)}(x)$ (dotted),
$\mathcal{K}_{2}^{(q)}(x)$ (dashed) and $\mathcal{K}_{0}^{(q)}(b)$ (solid
point) of the optimal strategy.%

\[%
\begin{array}
[c]{cc}%
{\parbox[b]{2.9362in}{\begin{center}
\includegraphics[
height=1.938in,
width=2.9362in
]%
{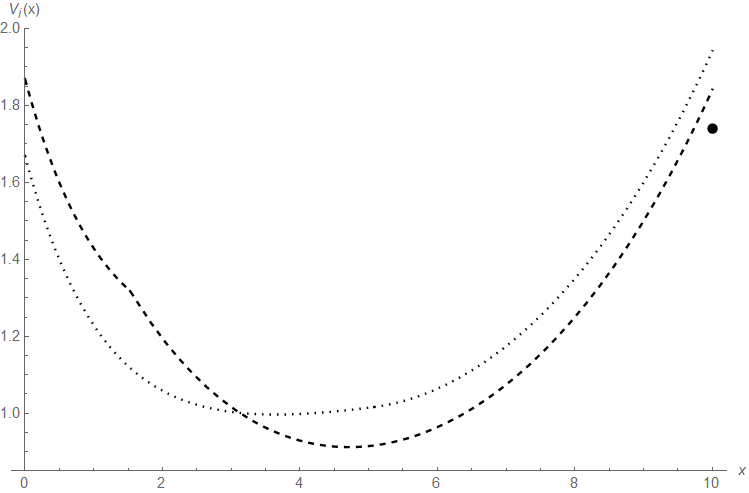}%
\\
Figure 2: First example, total cost.
\end{center}}}%
&
{\parbox[b]{2.9725in}{\begin{center}
\includegraphics[
height=1.9413in,
width=2.9725in
]%
{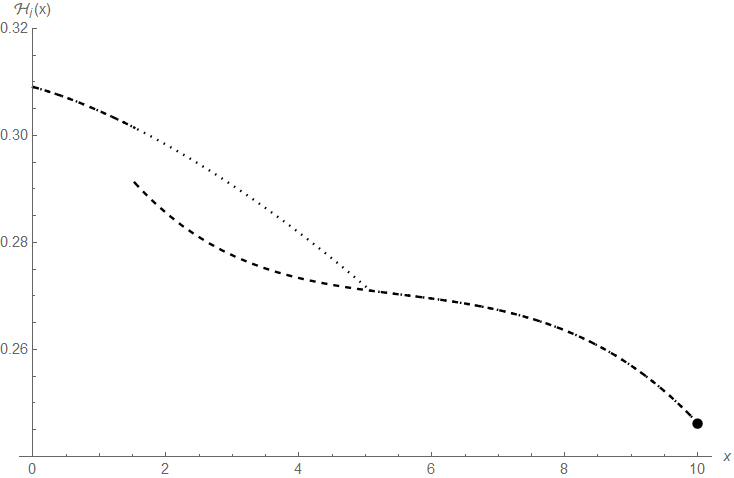}%
\\
Figure 3: First example, holding cost.
\end{center}}}%
\\%
{\parbox[b]{2.9346in}{\begin{center}
\includegraphics[
height=1.9339in,
width=2.9346in
]%
{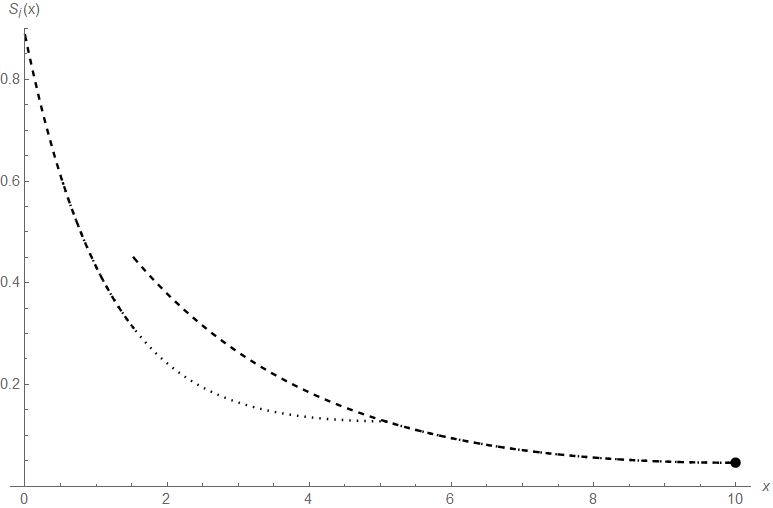}%
\\
Figure 4: First example, shortage cost.
\end{center}}}%
&
{\parbox[b]{2.8827in}{\begin{center}
\includegraphics[
height=1.9051in,
width=2.8827in
]%
{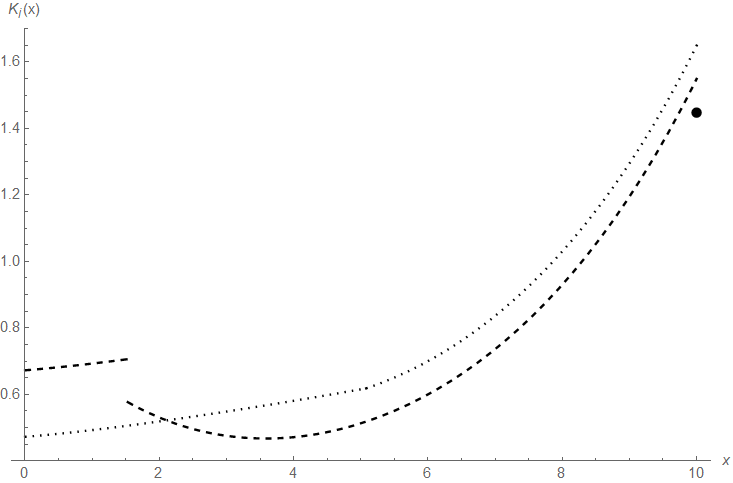}%
\\
Figure 5: First example, switching cost.
\end{center}}}%
\end{array}
\]

\begin{remark}
\label{Remark Examples}

(1) $V_{2}-V_{1}$ is equal to $K_{21}$ in $A_{{\small 21}}$ and equal to
$-K_{12}$ in $A_{{\small 12}}$; also $V_{2}(b^{-})-V_{0}(b)=K_{20}$ and
$V_{1}(b^{-})-V_{0}(b)=K_{12}+K_{10}$. $V_{2}$ is not differentiable at
$y_{{\small 2}},$ and so it is necessary to use the notion of viscosity solution.

(2) $\mathcal{H}_{2}^{(q)}=\mathcal{H}_{1}^{(q)}$ in the switching zones
$A_{{\small 21}}\cup$ $A_{{\small 12}}$. $\mathcal{H}_{2}^{(q)}$ is not
continuous at the boundary $y_{{\small 2}}$ between the switching zone
$A_{{\small 21}}$ and the non-action zone $(y_{{\small 2}},y_{{\small 1}}).$
The jump of $\mathcal{H}_{2}^{(q)}$ at $y_{{\small 2}}$ is downward because,
for an initial inventory level $x$ in the non-action zone, $X_{1,t}>$
$X_{2,t}$ for $t>0$ while these processes remain in the non-action zone. Also
note that $\mathcal{H}_{2}^{(q)}(b^{-})=\mathcal{H}_{1}^{(q)}(b^{-}%
)=\mathcal{H}_{0}^{(q)}(b).$

(3) $\mathcal{S}_{2}^{(q)}=\mathcal{S}_{1}^{(q)}$ in the switching zones
$A_{{\small 21}}\cup$ $A_{{\small 12}}$. As in the previous case and for
similar reasons, $\mathcal{S}_{2}^{(q)}$ has a discontinuity at $y_{{\small 2}%
}$, but in this case the jump is upward. Also note that $\mathcal{S}_{2}%
^{(q)}(b^{-})=\mathcal{S}_{1}^{(q)}(b^{-})=\mathcal{S}_{0}^{(q)}(b)$.

(4) $\mathcal{K}_{2}^{(q)}-\mathcal{K}_{1}^{(q)}$ is equal to $K_{21}$ in
$A_{{\small 21}}$ and equal to $-K_{12}$ in $A_{{\small 12}}$; also
$\mathcal{K}_{2}^{(q)}(b^{-})-\mathcal{K}_{0}^{(q)}(b)=K_{20}$ and
$\mathcal{K}_{1}^{(q)}(b^{-})-\mathcal{K}_{0}^{(q)}(b)=K_{12}+K_{20}$.
$\mathcal{K}_{2}^{(q)}$ has a downward jump at $y_{{\small 2}}$ because this
point is the boundary between the switching zone $A_{{\small 21}}$ and the
non-action zone $(y_{{\small 2}},y_{{\small 1}}).$
\end{remark}

\subsection{Second Example: Strategy of type one is optimal}

In this example, we consider that the demands are distributed as $Exp(1)$, the
parameters are $q=0.1$, $\lambda=2$, $l=0$, $b=20$, the rates of production
are $\sigma_{1}=2.5$, $\sigma_{2}=2.2$, and the costs are given by
$K_{12}=K_{21}=0.05$, $K_{20}=1/200$, $K_{10}=11/2000$,$~K_{01}=K_{02}=0$,
$h_{2}(x)=(20+x)/1000$, $h_{1}(x)=(30+x)/1000$, $h_{0}(b)=(2+10b)/10000$,
$p(y)=(80+40y)/100$.

In this case, the value functions of the best Doshi strategy do not satisfy
the condition of Theorem \ref{Teorema Verificacion}, so we look for the best
band strategies of type one, which is given by the sets $A_{{\small 12}%
}=[y_{{\small 1}},20)$, $A_{{\small 21}}=[0,y_{{\small 2}}]$, $C_{{\small 1}%
}=[0,y_{{\small 3}}]$ and $C_{{\small 2}}=(y_{{\small 3}},20)$ for
$y_{{\small 2}}=6.213$, $y_{3}=9.805$ and $y_{{\small 1}}=17.294$. The value
functions of this strategy of type one are viscosity solutions of the
equations (\ref{HJBi}) and satisfy the conditions of Theorem
\ref{Teorema Verificacion}, so this is the optimal strategy. We show this
optimal strategy in Figure 6.%

\[%
{\parbox[b]{4.8339in}{\begin{center}
\includegraphics[
height=2.6274in,
width=4.8339in
]%
{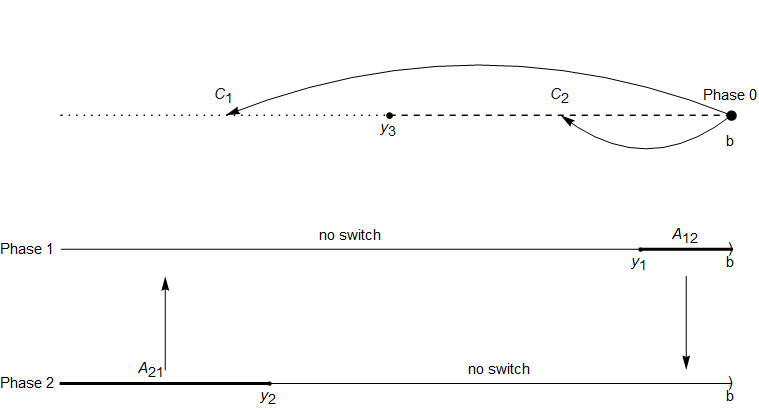}%
\\
Figure 6: Optimal strategy in second example.
\end{center}}}%
\]

In Figure 7, we show the discounted total cost $V_{1}(x)$ (dotted), $V_{2}(x)$
(dashed) and $V_{0}(b)$ (solid point) of the optimal strategy; in Figure 8,
the discounted holding cost $\mathcal{H}_{1}^{(q)}(x)$ (dotted),
$\mathcal{H}_{2}^{(q)}(x)$ (dashed) and $\mathcal{H}_{0}^{(q)}(b)$ (solid
point) of the optimal strategy; in Figure 9, the discounted penalty cost
$\mathcal{S}_{1}^{(q)}(x)$ (dotted), $\mathcal{S}_{2}^{(q)}(x)$ (dashed) and
$\mathcal{S}_{0}^{(q)}(b)$ (solid point) of the optimal strategy; and finally
in Figure 10, the discounted penalty cost $\mathcal{K}_{1}^{(q)}(x)$ (dotted),
$\mathcal{K}_{2}^{(q)}(x)$ (dashed) and $\mathcal{K}_{0}^{(q)}(b)$ (solid
point) of the optimal strategy.

The observations of Remark \ref{Remark Examples} hold for this example.
\[%
\begin{array}
[c]{cc}%
{\parbox[b]{2.8127in}{\begin{center}
\includegraphics[
height=1.8548in,
width=2.8127in
]%
{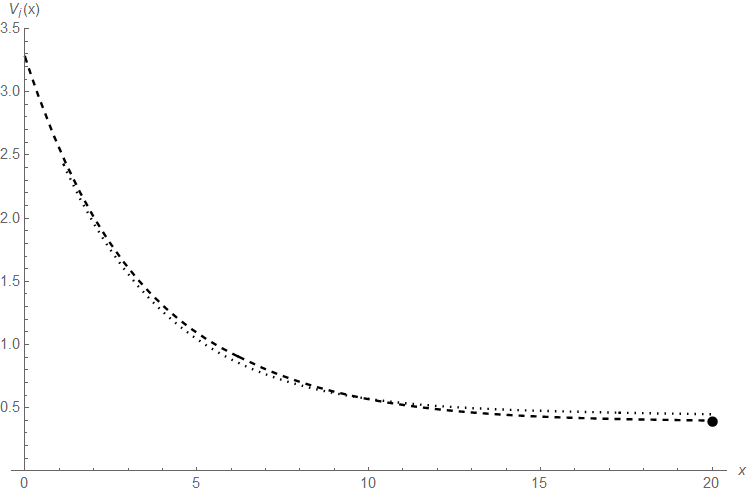}%
\\
Figure 7: Second example, total cost.
\end{center}}}%
&
{\parbox[b]{2.7979in}{\begin{center}
\includegraphics[
height=1.826in,
width=2.7979in
]%
{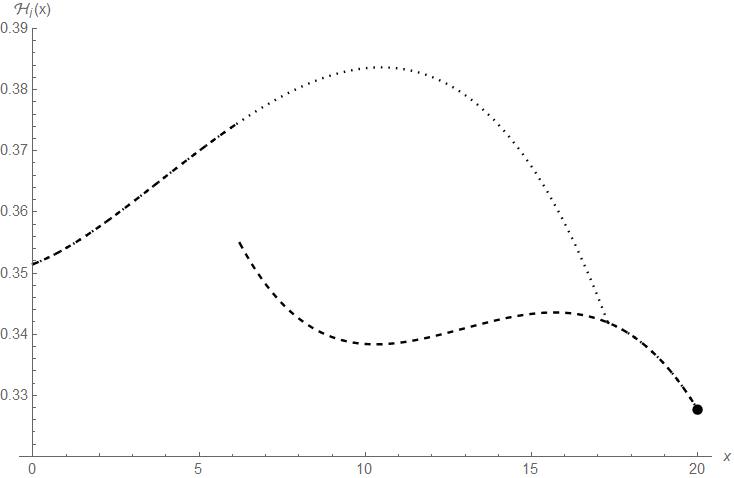}%
\\
Figure 8: Second example, holding cost.
\end{center}}}%
\\%
{\parbox[b]{2.8415in}{\begin{center}
\includegraphics[
height=1.8713in,
width=2.8415in
]%
{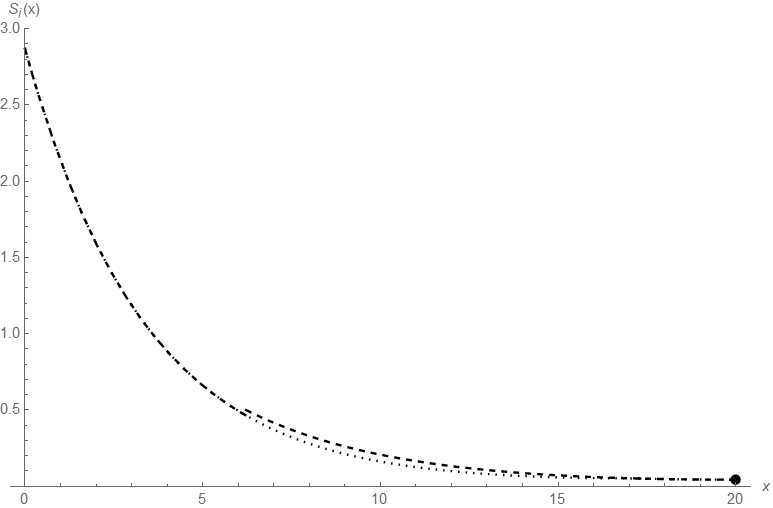}%
\\
Figure 9: Second example, shortage cost.
\end{center}}}%
&
{\parbox[b]{2.8045in}{\begin{center}
\includegraphics[
height=1.8482in,
width=2.8045in
]%
{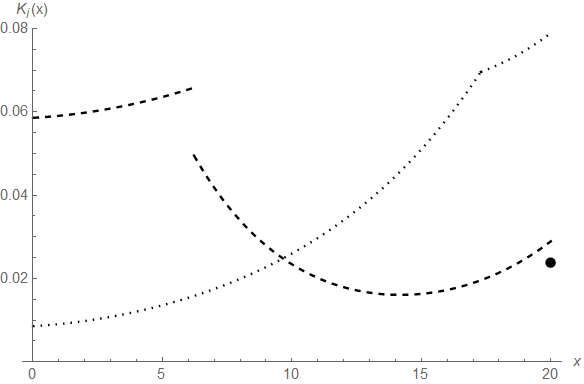}%
\\
Figure 10: Second example, switching cost.
\end{center}}}%
\end{array}
\]

\subsection{Third Example: Strategy of type two is optimal}

In this last example, we consider that the demands are distributed as
$Exp(1)$, the parameters are $q=0.1$, $\lambda=2$, $l=0$, $b=10$, the rates of
production are $\sigma_{1}=3.5$, $\sigma_{2}=2.5$, and the costs are given by
$K_{12}=K_{21}=0.05$, $K_{20}=0$, $K_{10}=1/100$,$~K_{01}=K_{02}=0$,
$h_{1}(x)=h_{2}(x)=(1+12x)/100$, $h_{0}(b)=(1+10b)/100$, $p(y)=2+1.1y$.

In this case, the value functions of the best strategy of type one do not
satisfy the condition of Theorem \ref{Teorema Verificacion}, so we look for
the best band strategies of type two, which is given by the sets
$A_{{\small 12}}=[y_{{\small 1}},y_{4}]$,\ $A_{{\small 21}}=[0,y_{{\small 2}%
}]$, $C_{{\small 1}}=[0,y_{{\small 3}}]$ and $C_{{\small 2}}=(y_{{\small 3}%
},10)$ for $y_{{\small 2}}=2.468$, $y_{3}=3.114$, $y_{{\small 1}}=4.610$,
$y_{{\small 4}}=7.660$. The value functions of this strategy of type two are
viscosity solutions of the equations (\ref{HJBi}) and satisfy the conditions
of Theorem \ref{Teorema Verificacion}, so this is the optimal strategy. We
show this optimal strategy in Figure 11.%

\[%
{\parbox[b]{4.8676in}{\begin{center}
\includegraphics[
height=2.5368in,
width=4.8676in
]%
{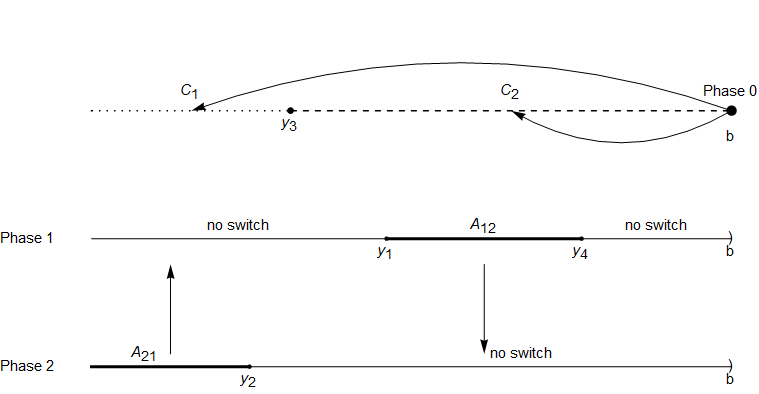}%
\\
Figure 11: Optimal strategy in third example.
\end{center}}}%
\]

In Figure 12, we show the discounted total cost $V_{1}(x)$ (dotted),
$V_{2}(x)$ (dashed) and $V_{0}(b)$ (solid point) of the optimal strategy; in
Figure 13, the discounted holding cost $\mathcal{H}_{1}^{(q)}(x)$ (dotted),
$\mathcal{H}_{2}^{(q)}(x)$ (dashed) and $\mathcal{H}_{0}^{(q)}(b)$ (solid
point) of the optimal strategy; in Figure 14, the discounted penalty cost
$\mathcal{S}_{1}^{(q)}(x)$ (dotted), $\mathcal{S}_{2}^{(q)}(x)$ (dashed) and
$\mathcal{S}_{0}^{(q)}(b)$ (solid point) of the optimal strategy; and finally
in Figure 15, the discounted penalty cost $\mathcal{K}_{1}^{(q)}(x)$ (dotted),
$\mathcal{K}_{2}^{(q)}(x)$ (dashed) and $\mathcal{K}_{0}^{(q)}(b)$ (solid
point) of the optimal strategy.

The observations (1), (2) and (3) of Remark \ref{Remark Examples} also hold
for this example. In this case $V_{1}$ is not differentiable at$y_{{\small 4}%
}.$ Also note, that $\mathcal{K}_{2}^{(q)}-\mathcal{K}_{1}^{(q)}$ is equal to
$K_{21}$ in $A_{{\small 21}}$ and equal to $-K_{12}$ in $A_{{\small 12}}$;
also $\mathcal{K}_{2}^{(q)}(b^{-})-\mathcal{K}_{0}^{(q)}(b)=K_{20}$ and
$\mathcal{K}_{1}^{(q)}(b^{-})-\mathcal{K}_{0}^{(q)}(b)=K_{10}$ because
$(y_{4},b)$ is the second component of the non-action zone. As in the previous
examples, $\mathcal{K}_{2}^{(q)}$ has a downward jump at $y_{{\small 2}}$ and,
in this case, $\mathcal{K}_{1}^{(q)}$ has a downward jump at $y_{{\small 4}%
}\ $because this point is the boundary between the switching zone
$A_{{\small 12}}=[y_{{\small 1}},y_{{\small 4}}]$ and the non-action zone
$(y_{{\small 4}},b).$
\[%
\begin{array}
[c]{cc}%
{\parbox[b]{2.9305in}{\begin{center}
\includegraphics[
height=1.9355in,
width=2.9305in
]%
{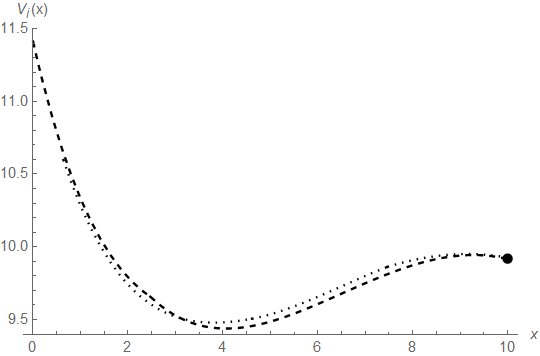}%
\\
Figure 12: Third example, total cost.
\end{center}}}%
&
{\parbox[b]{2.76in}{\begin{center}
\includegraphics[
height=1.8474in,
width=2.76in
]%
{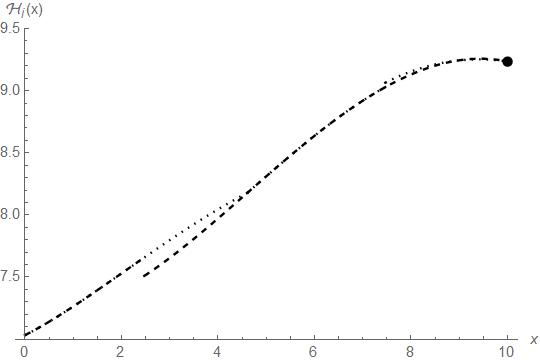}%
\\
Figure 13: Third example, holding cost.
\end{center}}}%
\\%
{\parbox[b]{2.788in}{\begin{center}
\includegraphics[
height=1.8935in,
width=2.788in
]%
{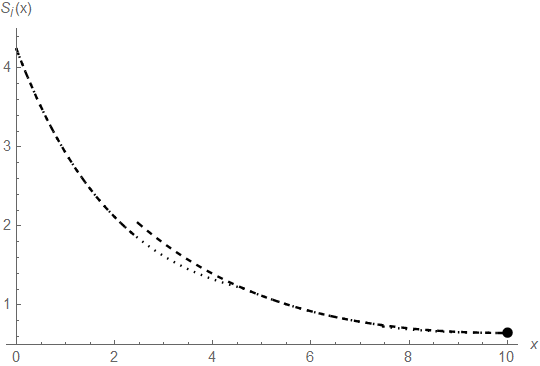}%
\\
Figure 14: Third example, shortage cost.
\end{center}}}%
&
{\parbox[b]{2.8292in}{\begin{center}
\includegraphics[
height=1.868in,
width=2.8292in
]%
{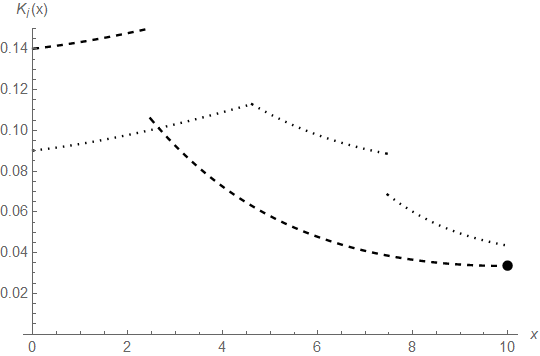}%
\\
Figure 15: Third example, switching cost.
\end{center}}}%
\end{array}
\]

\end{document}